\documentclass{amsart}
\usepackage{amscd,amssymb}
\usepackage{fullpage}
\usepackage{pb-diagram,lamsarrow,pb-lams}

\newcommand{\Aut}       {\operatorname{Aut}}
\newcommand{\Div}       {\operatorname{Div}}
\newcommand{\GL}        {\operatorname{GL}}
\newcommand{\Hom}       {\operatorname{Hom}}
\newcommand{\Level}     {\operatorname{Level}}
\newcommand{\Map}       {\operatorname{Map}}
\newcommand{\Nil}       {\operatorname{Nil}}
\newcommand{\Pic}       {\operatorname{Pic}}
\newcommand{\SDiv}      {\operatorname{SDiv}}
\newcommand{\Vect}      {\operatorname{Vect}}

\newcommand{\ind}       {\operatorname{ind}}
\newcommand{\res}       {\operatorname{res}}
\newcommand{\spec}      {\operatorname{spec}}
\newcommand{\spf}       {\operatorname{spf}}
\newcommand{\stab}      {\operatorname{stab}}
\newcommand{\tr}        {\operatorname{tr}}

\newcommand{\F}         {{\mathbb{F}}}
\newcommand{\N}         {{\mathbb{N}}}
\newcommand{\Z}         {{\mathbb{Z}}}
\newcommand{\Q}         {{\mathbb{Q}}}
\newcommand{\C}         {{\mathbb{C}}}
\newcommand{\Zp}        {{\mathbb{Z}_p}}          
\newcommand{\Zpl}       {{\mathbb{Z}_{(p)}}}      
\newcommand{\QZp}       {{\mathbb{Q}_p/\mathbb{Z}_p}}              

\newcommand{\al}        {\alpha}
\newcommand{\bt}        {\beta} 
\newcommand{\gm}        {\gamma}
\newcommand{\dl}        {\delta}
\newcommand{\ep}        {\epsilon}
\newcommand{\zt}        {\zeta}
\newcommand{\tht}       {\theta}

\newcommand{\sg}        {\sigma}
\newcommand{\kp}        {\kappa}
\newcommand{\lm}        {\lambda}
\newcommand{\om}        {\omega}
\newcommand{\tCh}       {\tau_{\text{Ch}}}
\newcommand{\rhoCh}     {\rho_{\text{Ch}}}
\newcommand{\xib}       {\overline{\xi}}

\newcommand{\Tht}       {\Theta}
\newcommand{\Lm}        {\Lambda}
\newcommand{\Sg}        {\Sigma}
\newcommand{\Om}        {\Omega}
\newcommand{\OmCh}      {\Omega_{\text{Ch}}}

\newcommand{\CA}        {{\mathcal{A}}}
\newcommand{\CC}        {{\mathcal{C}}}
\newcommand{\CF}        {{\mathcal{F}}}
\newcommand{\CG}        {{\mathcal{G}}}
\newcommand{\CGb}       {\overline{\mathcal{G}}}
\newcommand{\CGt}       {\widetilde{\mathcal{G}}}

\newcommand{\CN}        {{\mathcal{N}}}
\newcommand{\CNb}       {\overline{\mathcal{N}}}
\newcommand{\CNt}       {\widetilde{\mathcal{N}}}

\newcommand{\DD}        {{\mathbb{D}}}
\newcommand{\GG}        {{\mathbb{G}}}
\newcommand{\HH}        {{\mathbb{H}}}

\renewcommand{\O}       {\mathcal{O}}
\newcommand{\OG}        {\mathcal{O}_{\mathbb{G}}}
\newcommand{\OX}        {\mathcal{O}_X}
\newcommand{\OY}        {\mathcal{O}_Y}

\newcommand{\ot}        {\otimes}
\newcommand{\tm}        {\times}
\newcommand{\sse}       {\subseteq}
\newcommand{\st}        {\;|\;}
\newcommand{\sm}        {\setminus}

\newcommand{\bC}        {\overline{C}}
\newcommand{\bV}        {\overline{V}}
\newcommand{\xla}       {\xleftarrow}
\newcommand{\xra}       {\xrightarrow}
\newcommand{\Gmh}       {\widehat{\mathbb{G}}_m}
\newcommand{\XCh}       {X_{\text{Ch}}}
\newcommand{\cpi}       {{\mathbb{C}P^\infty}}
\newcommand{\hX}        {\widehat{\mathcal{X}}}
\newcommand{\haf}       {\widehat{\mathbb{A}}}
\newcommand{\hbT}       {\overline{h}\mathcal{T}}
\newcommand{\op}        {{\text{op}}}
\newcommand{\psb}[1]    {[\![#1]\!]}
\newcommand{\un}[1]     {\underline{#1}}
\newcommand{\ov}[1]     {\overline{#1}}
\newcommand{\aff}       {\mathbb{A}}
\newcommand{\convto}    {\Longrightarrow}
\newcommand{\hx}        {\hat{x}}
\newcommand{\ua}        {\underline{a}}
\newcommand{\ha}        {\hat{a}}
\newcommand{\ip}[2]     {\langle #1,#2\rangle}
\newcommand{\dts}       {\langle 2\rangle}
\newcommand{\pri}       {\mathfrak{p}}
\newcommand{\mxi}       {\mathfrak{m}}
\newcommand{\half}      {{\textstyle\frac{1}{2}}}
\newcommand{\hot}       {\widehat{\otimes}}
\newcommand{\ub}        {\underline{b}}
\newcommand{\unm}       {{\underline{m}}}
\newcommand{\unn}       {{\underline{n}}}
\newcommand{\cb}        {\overline{c}}
\newcommand{\xb}        {\overline{x}}
\newcommand{\yb}        {\overline{y}}
\newcommand{\zb}        {\overline{z}}

\newcommand{\colim}  {\operatornamewithlimits{\underset{\longrightarrow}{lim}}}

\newcommand{\bsm}       {\left(\begin{smallmatrix}}
\newcommand{\esm}       {\end{smallmatrix}\right)}
\newcommand{\bcf}[2]    {{\bsm #1\\#2\esm}}

\renewcommand{\:}{\colon}

\newtheorem{theorem}{Theorem}[section]

\newtheorem{lemma}[theorem]{Lemma}
\newtheorem{proposition}[theorem]{Proposition}
\newtheorem{corollary}[theorem]{Corollary}
\theoremstyle{definition}
\newtheorem{remark}[theorem]{Remark}
\newtheorem{definition}[theorem]{Definition}
\newtheorem{example}[theorem]{Example}

\newenvironment{diag}{
 \renewcommand{\typeout}[1]{}
 \begin{displaymath}
 \begin{diagram}}{
 \end{diagram}
 \end{displaymath}} 

\begin{document}
\title{Chern approximations for generalised group cohomology}
\author{N.~P.~Strickland}
\date{\today}
\bibliographystyle{abbrv}

\maketitle 

Let $G$ be a finite group, and let $E^*$ be a generalised cohomology
theory, subject to certain technical conditions (``admissibility'' in
the sense of~\cite{grst:vlc}).  Our aim in this paper is to define and
study a certain ring $C(E,G)$ that is in a precise sense the best
possible approximation to $E^0BG$ that can be built using only
knowledge of the complex representation theory of $G$.  There is a
natural map $C(E,G)\xra{}E^0BG$, whose image is the subring
$\bC(E,G)\leq E^0BG$ generated over $E^0$ by all Chern classes of all
complex representations.  There is ample precedent for considering
this subring in the parallel case of ordinary cohomology; see for
example~\cite{th:ccc,th:mrc,grle:scs}.  However, although the
generators of $\bC(E,G)$ come from representation theory, the same
cannot be said for the relations; one purpose of our construction is
to remedy this.  We also also develop a kind of generalised character
theory which gives good information about $\Q\ot C(E,G)$.  In the few
cases that we have been able to analyse completely, either $\Q\ot
C(E,G)\neq\Q\ot E^0BG$ for easy character-theoretic reasons, or we
have $C(E,G)=E^0BG$.

Rather than working directly with rings, we will study the formal
schemes $X(G)=\spf(E^0BG)$ and $\XCh(G)=\spf(C(E,G))$; note that the
map $C(E,G)\xra{}E^0BG$ corresponds to a map $X(G)\xra{}\XCh(G)$.
See~\cite{grst:vlc,st:fsf,st:fsfg} for foundational material on formal
schemes.  Suitably interpreted, our main definition is that $\XCh(G)$
is the scheme of homomorphisms from the $\Lm$-semiring $R^+(G)$ of
complex representations of $G$ to the $\Lm$-semiring scheme of
divisors on the formal group $\GG$ associated to $E$.

We start by fixing some conventions in Section~\ref{sec-notation}.  We
then recall the basic theory of $\Lm$-semirings
(Section~\ref{sec-LSR}), set up the parallel theory of $\Lm$-semiring
schemes, and define the $\Lm$-semiring scheme of divisors
(Section~\ref{sec-LSR-schemes}).  We then recall the definition of
Adams operations and study their basic properties
(Section~\ref{sec-Adams}).  Using this we give a precise definition of
$\XCh(G)$ and an implicit presentation of $C(E,G)$ by generators and
relations (Section~\ref{sec-chern-approx}).  In
Section~\ref{sec-Sigma-three}, we work out the case of the symmetric
group $\Sg_3$ at the prime $3$, and show that $X(\Sg_3)=\XCh(\Sg_3)$.
In Section~\ref{sec-Abelian} we show that $X(G)=\XCh(G)$ when $G$ is
Abelian, and in Section~\ref{sec-height-one} we show that the same is
true when $E$ is the $p$-adic completion of complex $K$-theory and $G$
is a $p$-group.  We then use Adams operations to reduce certain
questions to the Sylow subgroup of $G$ (Section~\ref{sec-Sylow}) and
to prove that $\XCh(G)$ is finite over $X=\spf(E^0)$
(Section~\ref{sec-finite}).  In Section~\ref{sec-gen-char}, we recall
the Hopkins-Kuhn-Ravenel generalised character theory, which relates
$\Q\ot E^0BG$ to the set $\Om(G)$ of conjugacy classes of
homomorphisms $\Z_p^n\xra{}G$.  We give a parallel (but less precise)
relationship between $\Q\ot C(E,G)$ and the set $\OmCh(G)$ of
$\Lm$-semiring homomorphisms $R^+(G)\xra{}\N[(\QZp)^n]$.  These
descriptions are related by a map $\kp\:\Om(G)\xra{}\OmCh(G)$.  In
Section~\ref{sec-OmChG}, we compare $\Om(G)$ and $\OmCh(G)$ with two
other sets that are sometimes easier to understand.  We next return to
examples: in Section~\ref{sec-Sigma-four} we show that
$X(\Sg_4)=\XCh(\Sg_4)$ at the prime $2$, and in
Section~\ref{sec-xspec} we study $\Om(G)$ and $\OmCh(G)$ when $G$ is
an extraspecial group at an odd prime.  We then show that a certain
approach which appears more precise actually captures no more
information (Section~\ref{sec-precise}).  We conclude in
Section~\ref{sec-restriction} by proving a result in representation
theory that was used in Section~\ref{sec-Sylow}.

\section{Notation and conventions}\label{sec-notation}

Fix a prime $p$.  Throughout this paper, $E$ will denote a $p$-local
generalised cohomology theory with an associative and unital product.
We write $E^k$ for $E^k(\text{point})$, so $E^*$ is a $\Z$-graded ring
and $E^0$ is an ungraded ring.  We assume that $E$ has the following
properties:
\begin{enumerate}
 \item $E^0$ is a commutative complete local Noetherian ring, with
  maximal ideal $\mxi$ say.
 \item $E^k=0$ whenever $k$ is odd.
 \item $E^{-2}$ contains a unit (so $E^kX\simeq E^{k-2}X$ for all
  $X$).
 \item Either $p>2$ and $E$ is commutative, or $p=2$ and $E$ is
  quasi-commutative, which means that there is a natural derivation
  $Q\:E^kX\xra{}E^{k+1}X$ and an element $v\in E^{-2}$ such that
  $2v=0$ and $ab-(-1)^{|a||b|}ba=vQ(a)Q(b)$ for all $a,b\in E^*X$.
\end{enumerate}
There is one more condition, which needs some background explanation.
Note that the quasi-commutativity condition means that whenever
$E^1X=0$, the ring $E^0X$ is commutative (in the usual ungraded
sense.)  In particular, $E^*=E^*(\text{point})$ is commutative.  A
collapsing Atiyah-Hirzebruch spectral sequence argument shows that
\[ E^*\cpi\simeq E^*\hot H^*\cpi=E^*\hot\Z\psb{y}=E^*\psb{y}, \]
where $y\in\widetilde{E}^2\cpi$; it follows that $E$ is
complex-orientable.  We can multiply $y$ by a unit in $E^{-2}$ to get
an element $x\in\widetilde{E}^0\cpi$ such that $E^*\cpi=E^*\psb{x}$.
We fix such an element $x$ once and for all (although we will state
our results in a form independent of this choice as far as possible).
This gives rise in the usual way to a formal group law $F$ over $E^0$.

\begin{enumerate}\setcounter{enumi}{4}
 \item In addition to the above properties, we will assume that
  the reduction of $F$ modulo the maximal ideal of $E^0$ has height
  $n<\infty$.  
\end{enumerate}

In the language of~\cite{grst:vlc}, our conditions say that
$E$ is a $K(n)$-local admissible theory.  The only difference is that
previously we insisted that $E$ should be commutative rather than
quasi-commutative; the reader can easily check that everything
in~\cite{grst:vlc} goes through in the quasi-commutative case.

We will describe our results in the language of formal schemes.  Most
of the formal schemes that we consider have the form $\spf(R)$, where
$R$ is a complete local Noetherian $E^0$-algebra.  For these the
foundational setting discussed in~\cite{st:fsf} is satisfactory: one
can regard the category of formal schemes as the opposite of the
category of complete semilocal Noetherian rings and continuous
homomorphisms.  We also make some use of formal schemes such as
$\spf(\prod_{k\in\Z}E^0\psb{c_1,c_2,\ldots})$; a set of foundations
covering these is developed in~\cite{st:fsfg}.  The older category of
formal schemes embeds as a full subcategory of the newer one.

\begin{definition}
 We let $X$ be the formal scheme $\spf(E^0)$, and write $\GG$ for the
 formal group $\spf(E^0\cpi)$ over $X$.  Note that our element
 $x\in\widetilde{E}^0\cpi$ can be regarded as a coordinate on $\GG$,
 with the property that
 \[ x(a+b) = F(x(a),x(b)) = x(a) +_F x(b). \]
\end{definition}

\begin{remark}
 Many of our constructions work with an arbitary formal group $\GG$
 over a formal scheme $X$; it is not usually necessary to assume that
 $\GG$ comes from a cohomology theory, although that is the case of
 most interest for us.
\end{remark}

We will let $G$ denote a finite group.  We write $e=e(G)$ for the
exponent of $G$, in other words the least common multiple of the
orders of the elements.  We factor this in the form
$e=p^ve'=p^{v(G)}e'(G)$, where $e'\neq 0\pmod{p}$.  We also choose a
Sylow $p$-subgroup $P\leq G$.

\section{$\Lm$-(semi)rings}\label{sec-LSR}

We will use the following definition:
\begin{definition}
 A \emph{semiring} is a set $R$ equipped with the following structure.
 \begin{itemize}
  \item A commutative and associative addition law with neutral
   element (written as $0$); we do not assume that there are additive
   inverses.
  \item A commutative and associative multiplication law with neutral
   element $1$, which distributes over addition.
 \end{itemize}

 A \emph{$\Lm$-semiring} is a semiring $R$ equipped with Maps
 $\lm^k\:R\xra{}R$ for $k\geq 0$ satisfying $\lm^0(x)=1$ and
 $\lm^1(x)=x$ and $\lm^k(x+y)=\sum_{k=i+j}\lm^i(x)\lm^j(y)$.

 A $\Lm$-ring is a $\Lm$-semiring which has additive inverses.  
\end{definition}

The initial $\Lm$-semiring is $\N$ and the initial $\Lm$-ring is $\Z$;
in both cases we have
\[ \lm^k(n)=\bsm n\\ k\esm= n(n-1)\ldots(n-k+1)/k!. \]

\begin{definition}
 An \emph{$\N$-augmented} $\Lm$-semiring is a $\Lm$-semiring $R$
 equipped with a homomorphism $\dim\:R\xra{}\N$ of $\Lm$-semirings.  A
 \emph{$\Z$-augmented} $\Lm$-ring is a $\Lm$-ring $R$ equipped with a
 homomorphism $\dim\:R\xra{}\Z$ of $\Lm$-rings.
\end{definition}

\begin{example}\label{eg-RG}
 Let $R^+(G)$ be the semiring of isomorphism classes of complex
 representations of $G$.  It is well-known that this is a
 $\Lm$-semiring with operations $\lm^k$ given by exterior powers.
 There is an augmentation $\dim\:R^+(G)\xra{}\N$ sending each
 representation to its dimension.
\end{example}

\begin{example}\label{eg-NA}
 Let $A$ be an Abelian group, and let $\N[A]$ be the group semiring of
 $A$, in other words the set of expressions $\sum_an_a[a]$ with
 $n_a\in\N$ and $n_a=0$ for all but finitely many $a$.  Equivalently,
 we have $\N[A]=\coprod_nA^n/\Sg_n$.  This has a canonical structure
 as a $\Lm$-semiring, with
 \[ \lm^k([a_1]+\ldots+[a_n]) = \sum_I [a_{i_1}+\ldots+a_{i_k}], \]
 where the sum on the right runs over all lists $I=(i_1,\ldots,i_k)$
 such that $1\leq i_1<\ldots<i_k\leq n$.  There is an augmentation
 $\dim\:\N[A]\xra{}\N$ defined by 
 \[ \dim([a_1]+\ldots+[a_n]) = n. \]
 If $A$ is finite and $A^*=\Hom(A,S^1)$ then $\N[A]=R^+(A^*)$ as
 $\Lm$-semirings.
\end{example}

\begin{example}\label{eg-Vect}
 For any space $X$, let $\Vect^+(X)$ denote the semiring of
 isomorphism classes of complex vector bundles over $X$.  This is a
 $\Lm$-semiring with operations as for $R^+(G)$.  We will always allow
 vector bundles to have different dimensions over different components
 of the base, so we do not have a natural map
 $\dim\:\Vect^+(X)\xra{}\Z$.  We write $\Vect_d^+(X)$ for the set of
 isomorphism classes of bundles all of whose fibres have dimension
 $d$, and we put $\Pic(X)=\Vect_1(X)\simeq H^2(X)$.  This is an
 Abelian group, and there is an evident map
 $\N[\Pic(X)]\xra{}\Vect^+(X)$.  In the case $X=BG$, there is a
 well-known homomorphism $R^+(G)\xra{}\Vect^+(BG)$ sending a
 representation $V$ to the bundle $V\tm_GEG$.
\end{example}

\begin{remark}\label{rem-identities}
 In the important examples of $\Lm$-(semi)rings, some extra identities
 hold that relate the elements $\lm^i\lm^j(x)$ and $\lm^i(xy)$ to the
 elements $\lm^k(x)$ and $\lm^l(y)$.  For many purposes it would be
 preferable to take these identities as part of the definition of a
 $\Lm$-(semi)ring.  However, it turns out that this would make no
 difference for us and the identities are complicated (particularly in
 the semiring case) so we omit them.  In Section~\ref{sec-precise} we
 will discuss an approach which is apparently even more precise, and
 show that it actually gives no more information than our approach
 using $\Lm$-semirings without extra identities.
\end{remark}

\begin{remark}\label{rem-Grothendieck}
 Let $R^+$ be a $\Lm$-semiring, and let $R$ be its Grothendieck
 completion, or in other words the group completion of $R^+$
 considered as a monoid under addition.  It is well-known that this
 can be made into a $\Lm$-ring in a canonical way, and that any
 homomorphism from $R^+$ to a $\Lm$-ring factors uniquely through $R$.
 Moreover, if $R^+$ is augmented over $\N$ then $R$ is augmented over
 $\Z$.
\end{remark}

\begin{example}\label{eg-Grothendieck}
 The Grothendieck completion of $R^+(G)$ is of course the ring $R(G)$
 of virtual representations of $G$, and the completion of $\N[A]$ is
 the group ring $\Z[A]$.  We write $\Vect(X)$ for the Grothendieck
 completion of $\Vect^+(X)$.  It is well-known that the complex
 $K$-theory $K^0(X)$ is a $\Z$-augmented $\Lm$-semiring and that there
 is a natural map $\Vect(X)\xra{}K^0(X)$ which is an isomorphism
 whenever $X$ is compact Hausdorff.
\end{example}

\begin{remark}
 We will occasionally use the notation $\Z[A]^+=\N[A]$ and
 $\Z[A]_d^+=A^d/\Sg_d\subset\N[A]$.
\end{remark}

\section{$\Lm$-(semi)ring schemes}\label{sec-LSR-schemes}

The theory of $\Lm$-semirings is an instance of universal algebra: it
is defined in terms of operations $\om\:R^k\xra{}R$ with $k=0,1$ or
$2$, and identities between operations derived from these.  It is thus
formal to define the notion of a $\Lm$-semiring object in any category
$\CC$ with finite products: such a thing is an object $R\in\CC$
equipped with maps 
\begin{align*}
 0,1   \: 1   & \xra{} R \\ 
 +,\tm \: R^2 & \xra{} R \\
 \lm^k \: R   & \xra{} R \hspace{6em} \text{ for all } k\in\N 
\end{align*}
making the evident diagrams commute.  (Here the object $1\in\CC$ is
the terminal object.)  Similar remarks apply to $\Lm$-rings.

Next, suppose that $\CC$ has arbitrary coproducts such that the
natural map 
\[ \coprod_{i,j}X_i\tm Y_j\xra{}\coprod_iX_i\tm\coprod_jX_j \]
is always an isomorphism.  We then have a product-preserving functor
$S\mapsto\un{S}:=\coprod_{s\in S}1$ from sets to $\CC$, so $\un{\N}$
is a $\Lm$-semiring object in $\CC$.  Similarly, $\un{\Z}$ is a
$\Lm$-ring object in $\CC$.  

\begin{example}\label{eg-BU}
 Take $\CC=\hbT$, the homotopy category of unbased CW-complexes.  We
 have a functor $\Vect^+(-)$ from $\hbT^\op$ to the category of
 $\Lm$-semirings, which is represented by the space $\coprod_dBU(d)$.
 It follows by Yoneda's lemma that $\coprod_dBU(d)$ is a
 $\Lm$-semiring in $\hbT$.  Similarly, the functor $K^0(-)$ from
 $\hbT^\op$ to the category of $\Lm$-rings is represented by the
 $\Lm$-ring space $\Z\tm BU$.  Note that in this context the object
 $\un{\N}$ is just the discrete space $\N$ and similarly for
 $\un{\Z}$, so $\coprod_dBU(d)$ is augmented over $\un{\N}$ and 
 $\Z\tm BU$ is augmented over $\un{\Z}$.
\end{example}

Now let $X$ be a formal scheme, and consider the category $\hX_X$ of
category of formal schemes over $X$ in the sense of~\cite{st:fsfg}.
For simplicity we will assume that $X$ is solid, which means that
$X=\spf(\OX)$ for some formal ring $\OX$.  Let $\CA$ be the category
of discrete $\OX$-algebras, and let $\CF$ be the category of functors
from $\CA$ to sets.  The category $\hX_X$ can be regarded as a
subcategory of $\CF$ (compare~\cite[Remark~2.1.5]{st:fsfg}), and the
inclusion $\hX_X\xra{}\CF$ preserves products.

We will refer to $\Lm$-(semi)ring objects in $\hX_X$ as
$\Lm$-(semi)ring schemes (suppressing the words ``formal'' and ``over
$X$'' for brevity).

Let $\GG$ be an ordinary formal group over $X$, in other words a
commutative group object in $\hX_X$ that is isomorphic in $\hX_X$ to
$\haf^1_X=\spf(\OX\psb{x})$.  We can then define the schemes
\begin{align*}
 \Div_d^+(\GG) &= \GG^d/\Sg_d                           \\
 \Div^+(\GG)   &= \coprod_{d\in\N} \Div_d^+(\GG)        \\
 \Div_0(\GG)   &= \colim_d\Div_d^+(\GG)                 \\
 \Div(\GG)     &= \un{\Z}\tm\Div_0(\GG)                 \\
 \Div_d(\GG)   &= \{d\}\tm\Div_0(\GG)\subset\Div(\GG).
\end{align*}
More detailed definitions are given in~\cite[Section 5]{st:fsfg},
where it is also explained how these formal schemes relate to the
theory of divisors on $\GG$.  In~\cite[Proposition 6.2.7]{st:fsfg} it
is observed that 
\begin{enumerate}
 \item $\Div^+(\GG)$ is the free commutative monoid object in $\CC$
  generated by $\GG$.
 \item $\Div(\GG)$ is the free commutative group object in $\CC$
  generated by $\GG$.
 \item $\Div_0(\GG)$ is the free commutative monoid object generated
  by $\GG$ considered as a based object in $\CC$, which is the same as
  the free commutative group object generated by $\GG$ considered as a
  based object in $\CC$.
\end{enumerate}
Moreover, all these universal properties are stable under base change:
if $X'$ is a formal scheme over $X$ then $\Div^+(\GG)\tm_XX'$ is the
free commutative monoid in $\hX_{X'}$ generated by $\GG\tm_XX'$ and so
on.  

Recall that $\OG=\OX\psb{x}$ and thus
$\O_{\GG^d}=\OX\psb{x_1,\ldots,x_d}$.  If $c_k$ denotes the
coefficient of $t^{d-k}$ in $\prod_i(t-x_i)$ then
$\O_{\Div^+_d(\GG)}=\OX\psb{c_1,\ldots,c_d}$ and
$\O_{\Div^+(\GG)}=\prod_{d\geq 0}\OX\psb{c_1,\ldots,c_d}$.  There are
also isomorphisms
\begin{align*}
 \O_{\Div_0(\GG)} &= \OX\psb{c_1,c_2,\ldots} \\
 \O_{\Div(\GG)}   &= \prod_{d\in\Z}\OX\psb{c_1,c_2,\ldots}.
\end{align*}
Using these, one sees that $\Div^+_d(\GG)$ is a closed subscheme of
$\Div_d(\GG)$, and $\Div^+(\GG)$ is a closed subscheme of
$\Div(\GG)$. 

If $E$ is an even periodic ring spectrum, $X=\spec(\pi_0E)$ and
$\GG=\spf(E^0\cpi)$ then there are natural isomorphisms
\begin{align*}
 \spf(E^0BU(d))      &= \Div_d^+(\GG) \\
 \spf(E^0BU)         &= \Div_0(\GG) \\
 \spf(E^0(\Z\tm BU)) &= \Div(\GG).
\end{align*}
This is just a translation of well-known calculations; details are
given in~\cite[Section 8]{st:fsfg}.

\begin{proposition}\label{prop-Div}
 Let $\GG$ be an ordinary formal group over a scheme $X$.  Then
 $\Div^+(\GG)$ has a natural structure as a $\Lm$-semiring scheme, and
 $\Div(\GG)$  has a natural structure as a $\Lm$-ring scheme.
 Moreover, there is a canonical homomorphism
 $\dim\:\Div(\GG)\xra{}\un{\Z}$ of $\Lm$-ring schemes, which sends
 $\Div_d(\GG)$ to $d$.
\end{proposition}
\begin{proof}
 Recall that $\CF$ is the category of functors from discrete
 $\OX$-algebras to sets.  Define $R^+,R\in\CF$ by $R^+(A)=\N[\GG(A)]$
 and $R(A)=\Z[\GG(A)]$.  It is clear that $R^+$ is a $\Lm$-semiring
 object in $\CF$, and $R$ is a $\Lm$-ring object.

 There is an evident inclusion
 $j\:\GG=\Div^+_1(\GG)\xra{}\Div^+(\GG)$.  As $\Div^+(\GG)$ is a
 commutative monoid scheme, the set $\Div^+(\GG)(A)$ is a commutative
 monoid for all $A\in\CA$.  As $R^+(A)$ is the free commutative monoid
 generated by the set $\GG(A)$, there is a unique homomorphism
 $\phi^+\:R^+(A)\xra{}\Div^+(\GG)(A)$ extending $j$.  These maps are
 natural in $A$ so we get a map $\phi^+\:R^+\xra{}\Div^+(\GG)$ in
 $\CF$.  If we interpret the colimits in $\hX_X$ then we have
 $\Div^+(\GG)=\coprod_d\GG^d/\Sg_d$; this translates to the statement
 that $\Div^+(\GG)$ is the initial example of a formal scheme
 in $\hX_X$ equipped with a map $R^+\xra{}\Div^+(\GG)$ in $\CF$.  By
 similar arguments, we find that $\Div(\GG)$ is the initial example of
 a formal scheme over $X$ with a map $\phi\:R\xra{}\Div(\GG)$ in
 $\CF$.  Moreover, one can check that the schemes $\Div^+(\GG)^k$ and
 $\Div(\GG)^k$ enjoy the evident analogous universal properties for
 all $k\geq 0$.

 It now follows that there is a unique map
 $\tm\:\Div^+(\GG)\tm\Div^+(\GG)\xra{}\Div^+(\GG)$ making the
 following diagram commute:
 \begin{diag}
  \node{R^+\tm R^+} \arrow{s,l}{\phi^+\tm\phi^+} \arrow{e,t}{\tm} 
  \node{R^+}        \arrow{s,r}{\phi^+}                           \\
  \node{\Div^+(\GG)\tm\Div^+(\GG)}               \arrow{e,b}{\tm}
  \node{\Div^+(\GG).}
 \end{diag}
 Similarly, all the other structure maps for the $\Lm$-semiring
 structure on $R^+$ induce operations on $\Div^+(\GG)$, and one checks
 easily that this makes $\Div^+(\GG)$ into a $\Lm$-semiring scheme.  A
 similar argument works for $\Div(\GG)$.  It is clear that there is a
 map $\dim\:\Div(\GG)\xra{}\un{\Z}$ as described.
\end{proof}

The above $\Lm$-semiring structure can be made more explicit as
follows.  Let $c_{d,k}\in\OX\psb{x_1,\ldots,x_d}$ be defined by
\[ \prod_{i=1}^d(t-x_i) = \sum_{i=0}^d c_{d,i} x^{d-i}. \]
Let $p_{d,e,k}(c_{d,1},\ldots,c_{d,d},c'_{e,1},\ldots,c'_{e,e})$ be
defined by
\[ \prod_{i=1}^d \prod_{j=1}^e (t-(x_i+_Fx'_j)) =
    \sum_{k=0}^{de} p_{d,e,k} t^{de-k}.
\]
Suppose $d,r\in\N$ and put $N=\bcf{d}{r}$.  Let
$q_{d,r,k}(c_{N,1},\ldots,c_{N,N})$ be defined by
\[ \prod_I(t-\sum^F_jx_{i_j}) = \sum_{k=0}^N q_{d,r,k} t^{N-k}, \]
where the sum on the left runs over all lists $I=(i_1,\ldots,i_r)$
such that $1\leq i_1<\ldots<i_r\leq d$.  Then the multiplication map
\[ \tm\:\Div_d^+(\GG)\tm\Div_e^+(\GG)\xra{}\Div_{de}^+(\GG) \]
corresponds to the map 
\[ \OX\psb{c_{de,1},\ldots,c_{de,de}} \xra{}
   \OX\psb{c_{d,1},\ldots,c_{d,d},c'_{1,e},\ldots,c'_{e,e}} 
\]
(of formal $\OX$-algebras) sending $c_{de,k}$ to $p_{d,e,k}$.
Similarly, the map corresponding to
$\lm^r\:\Div_d^+(\GG)\xra{}\Div_N^+(\GG)$ sends $c_{N,k}$ to
$q_{d,r,k}$.

\section{Adams operations}\label{sec-Adams}

We now recall the theory of Adams operations in $\Lm$-semirings; for a
more detailed exposition see~\cite{hu:fb}, for example.

Let $R$ be a $\Lm$-ring.  For any $a\in R$ we can form the power
series
\[ \lm_t(a) = \sum_{k\geq 0}\lm^k(a) (-t)^k \in R\psb{t}. \]
This is equal to $1$ mod $t$ and thus is invertible in $R\psb{t}$.  It
is easy to check that $\lm_t(0)=1$ and
$\lm_t(a+b)=\lm_t(a)\lm_t(b)$.  

We next define 
\[ \psi_t(a) = -t \lm_{-t}(a)^{-1}\, d\lm_{-t}(a)/dt \in R\psb{t}, \]
and let $\psi^k(a)$ be the coefficient of $t^k$ in $\psi_t(a)$.  This 
defines an additive map $\psi^k\:R\xra{}R$, called the $k$'th Adams
operation. 

Now consider the case $R=\Z[A]$ for some Abelian group $A$.  It is not
hard to see that 
\[ \psi^k(\sum_i n_i[a_i]) = \sum_i n_i[a_i]^k = \sum_i n_i[ka_i], \]
so $\psi^k$ is just the map induced by the homomorphism
$k.1_A\:A\xra{}A$.  Thus, if $A$ is actually a $\Zpl$-module or a
$\Zp$-module, then there is a natural way to define
$\psi^k\:\Z[A]\xra{}\Z[A]$ for all $k\in\Zpl$ or $k\in\Zp$ as
appropriate.  Moreover, we see that $\psi^k$ is a ring homomorphism
which preserves the semiring $\Z[A]^+$ and the subsets $\Z[A]_d^+$,
and that $\psi^k\psi^j=\psi^{kj}$ and $\psi^k\lm^j=\lm^j\psi^k$.

Now consider the $\Lm$-ring scheme $\Div(\GG)$.  As our original
definition of $\psi^k$ is natural, we evidently get morphisms 
\[ \psi^k\:\Div(\GG) \xra{} \Div(\GG) \]
of schemes.  It is well-known that $\GG$ is actually a $\Zp$-module
scheme, or in terms of our coordinate, that one can define the series
$[k]_F(x)$ in a sensible way for all $k\in\Zp$.  This means that each
ring $\Z[\GG(A)]$ admits Adams operations $\psi^k$ for all $k\in\Zp$,
with properties as above.  The argument of Proposition~\ref{prop-Div}
shows that
\begin{itemize}
 \item We can define operations $\psi^k$ on $\Div(\GG)$ for all
  $k\in\Zp$, extending the definition given previously.
 \item These maps are maps of ring schemes, induced by the maps
  $k\:\GG\xra{}\GG$.
 \item We have $\psi^j\psi^k=\psi^{jk}$ for all $j,k\in\Zp$, and
  $\psi^k\lm^j=\lm^j\psi^k$ for all $k\in\Zp$ and $j\in\N$.
 \item The map $\psi^k$ preserves $\Div_d(\GG)$, $\Div^+(\GG)$ and
  $\Div_d^+(\GG)=\GG^d/\Sg_d$ for all $d$.
\end{itemize}
 
\begin{lemma}\label{lem-p-torsion}
 For any discrete $\OX$-algebra $A$, the group $\GG(A)$ is a
 $p$-torsion group.
\end{lemma}
\begin{proof}
 Our coordinate $x$ gives an isomorphism
 $x\:\GG(A)\xra{}\haf^1(A)=\Nil(A)$ (the set of nilpotents in $A$).
 As $p$ lies in the maximal ideal of $E^0=\OX$ and $A$ is a discrete
 $\OX$-algebra, we see that $p^r=0\in A$ for some $r$, and thus
 $[p^r](x)$ is divisible by $x^2$ in $A\psb{x}$.  It follows that
 $[p^{rs}](x)$ is divisible by $x^{2^s}$.  For any $a\in\GG(A)$ we
 have $x(a)^{2^s}=0$ for large $s$, so $x(p^{rs}a)=0$ for large $s$,
 so $p^ma=0$ for large $m$ as required.
\end{proof}

\begin{lemma}\label{lem-psi-torsion}
 Let $A$ be a discrete $\OX$-algebra, and suppose we have a divisor 
 $D\in\Div(\GG)(A)$.  Then $\psi^k(D)=\dim(D)[0]$
 whenever the $p$-adic valuation $v_p(k)$ is sufficiently large.
\end{lemma}
\begin{proof}
 First suppose that $D\in\Div_d^+(\GG)(A)$.  We can then choose a
 faithfully flat map $A\xra{}A'$ such that the image of $D$ in
 $\Div_d^+(\GG)(A')$ has the form $\sum_{i=1}^d[a_i]$.  The map
 $\Div_d^+(\GG)(A)\xra{}\Div_d^+(\GG)(A')$ is automatically injective,
 so it suffices to show that for large $m$ we have
 $\psi^{p^m}(\sum_i[a_i])=\sum_i[p^ma_i]=d[0]$, which is immediate
 from the previous lemma.  

 Now suppose that $D\in\Div_d(\GG)(A)$.  We can then write $D$ in the
 form $D'-e[0]$ for some $D'\in\Div_{d+e}^+(\GG)(A)$ and we reduce
 easily to the previous case.

 Finally, consider a general divisor $D\in\Div_d(\GG)(A)$, which need
 not have constant dimension.  Instead, we have a splitting
 $A=A_1\tm\ldots\tm A_r$ giving a bijection
 $\Div(\GG)(A)=\prod_i\Div(\GG)(A_i)$ under which $D$ becomes an
 $r$-tuple $(D_1,\ldots,D_r)$ with $D_i\in\Div_{d_i}^+(\GG)(A_i)$ for
 some integers $d_i$.  This means that $\dim(D)$ becomes
 $(d_1,\ldots,d_r)$ under the bijection
 $\un{\Z}(A)=\prod_i\un{\Z}(A_i)$.  The cases considered previously
 imply that 
 \[ \psi^kD = (\psi^kD_1,\ldots,\psi^kD_r) = 
    (d_1[0],\ldots,d_r[0]) = \dim(D)[0]
 \]
 when $v_p(k)\gg 0$, as required.
\end{proof}

Now consider instead the $\Lm$-ring $R(G)$.  In this case we can
define Adams operations $\psi^k$ for $k\in\N$, and it is well-known
that in terms of characters we have
\[ \chi_{\psi^kV}(g) = \chi_V(g^k).  \]
As a virtual representation is determined by its character and
$\dim(V)=\chi_V(1)$, it follows easily that $\psi^k$ is a
degree-preserving map of $\Lm$-rings and that
$\psi^j\psi^k=\psi^{jk}$.  Moreover, if $e$ is the exponent of $G$ (in
other words, least common multiple of the orders of the elements) then
$\psi^k$ depends only on the congruence class of $k$ modulo $e$.  If
$k$ is coprime to $e$ (or equivalently, to $|G|$), it follows that
$\psi^k\psi^j=1$ for some $j$, so $\psi^k$ is an isomorphism.  In this
case the map $g\mapsto g^k$ is a bijection, and it follows easily that
$\psi^k$ preserves the usual inner product on $R(G)$.  A virtual
representation $V$ is an irreducible honest representation iff
$\chi_V(1)>0$ and $\ip{V}{V}=1$, and it follows that $\psi^k$ sends
irreducibles to irreducibles and thus sends $R^+_d(G)$ to $R^+_d(G)$.
(Compare~\cite[Exercise 9.4]{se:lrf}.)

However, if $k$ is not coprime to $e$ then $\psi^k$ need not preserve
$R^+(G)$.  For example, take $G=\Sg_3$, let $\ep$ be the nontrivial
one-dimensional representation, and let $\rho$ be the irreducible
two-dimensional representation.  We then have
$\psi^2(\rho)=\rho+1-\ep\not\in R^+(G)$.

\begin{lemma}\label{lem-psipv}
 Let $p^v$ be the $p$-part of the exponent of $G$.  Then for any
 homomorphism $f\:R(G)\xra{}\Div(\GG)(A)$ of $\Lm$-rings and any
 $V\in R_d(G)$ we have $f(V)\in\Div_d(\GG)$ and
 $\psi^{p^v}f(V)=d[0]$.
\end{lemma}
\begin{proof}
 Let the exponent of $G$ be $e=p^ve'$, where $e'$ is coprime to $p$.
 The map $\psi^{e'}\:\Div(\GG)\xra{}\Div(\GG)$ is an isomorphism and
 fixes $d[0]$, and $\psi^{e'}\psi^{p^v}=\psi^e$ so it suffices to show
 that $\psi^ef(V)=d[0]$.  To see this note that
 $\psi^ef(V)=f(\psi^eV)$ and
 $\chi_{\psi^e V}(g)=\chi_V(g^e)=\chi_V(1)=d$ for all $g$, so
 $\psi^eV$ is the trivial representation of rank $d$.  As $f$ is a
 ring map, we have $f(\psi^eV)=f(d)=d[0]$, as required.

 This implies that $\psi^{p^k}f(V)=d[0]$ for $k\gg 0$ but
 Lemma~\ref{lem-psi-torsion} says that $\psi^{p^k}f(V)=\dim(f(V))[0]$
 for $k\gg 0$, so $\dim(f(V))=d$, so $f(V)\in\Div_d(\GG)$ as claimed.
\end{proof}

\section{Chern approximations}\label{sec-chern-approx}

\begin{definition}
 Let $G$ be a finite group, and let $A$ be a discrete $\OX$-algebra.
 We define a functor $\XCh(G)$ from discrete $\OX$-algebras to sets by
 \[ \XCh(G)(A) =
     \{\text{ homomorphisms $R^+(G)\xra{}\Div^+(\GG)(A)$ 
              of $\Lm$-seimirings }\}.
 \]
 We write $C(E,G)$ for the ring $\O_{\XCh(G)}$ of natural
 transformations from $\XCh(G)$ to the forgetful functor $\aff^1$.  We
 also put $X(G)=\spf(E^0BG)$.  We refer to $C(E,G)$ as the Chern
 approximation to $E^0BG$, and to $\XCh(G)$ as the Chern approximation
 to $X(G)$.
\end{definition}

\begin{remark}
 We say that a homomorphism $f\:R(G)\xra{}\Div(\GG)(A)$ of $\Lm$-rings
 is \emph{positive} if $f(R^+(G))\sse\Div^+(\GG)(A)$.  It is clear
 from Remark~\ref{rem-Grothendieck} that $\XCh(G)(A)$ bijects
 naturally with the set of positive homomorphisms
 $R(G)\xra{}\Div(\GG)(A)$, and we will implicitly use this
 identification where convenient.  We also see from
 Lemma~\ref{lem-psipv} that positive homomorphisms satisfy
 $f(R^+_d(G))\sse\Div_d^+(\GG)(A)$.
\end{remark}

\begin{proposition}\label{prop-scheme}
 The functor $\XCh(G)$ is a formal scheme over $X$.  The ring
 $C(E,G)=\O_{\XCh(G)}$ is a quotient of a formal power series ring in
 finitely many variables over $\OX$ (and thus is a complete Noetherian
 local ring).
\end{proposition}
\begin{proof}
 Let $V_1,\ldots,V_h$ be the irreducible representations of $G$, and
 let $d_1,\ldots,d_h$ be their degrees.  We assume that these are
 ordered so that $V_1$ is the trivial representation of rank one.
 There are then natural numbers $m_{ijk}$ and $l^r_{ij}$ for $r\geq 0$
 and $1\leq i,j,k\leq h$ such that 
 \begin{align*}
  V_i\ot V_j &\simeq \bigoplus_k m_{ijk}.V_k    \\
  \lm^r V_i  &\simeq \bigoplus_j l^r_{ij}. V_j
 \end{align*}
 (Here $m.W$ means the direct sum of $m$ copies of $W$.)

 To give a homomorphism $f\:R^+(G)\xra{}\Div^+(\GG)(A)$ is the same as
 to give divisors $D_i=f(V_i)\in\Div_{d_i}^+(\GG)$ for $i=1,\ldots,h$
 such that
 \begin{align*}
  D_iD_j   &= \sum_{ijk} m_{ijk}D_k \\
  \lm^rD_i &= \sum_j l^r_{ij} D_j
 \end{align*}
 This exhibits $\XCh(G)(A)$ as the equaliser of a pair of maps from
 $\prod_{i=1}^h\Div_{d_i}(\GG)$ to
 \[ \prod_{i,j}\Div_{d_id_j}(\GG) \tm
    \prod_{r,i}\Div_{\bcf{d_i}{r}}(\GG).
 \]
 In particular, this is a pair of maps between formal schemes over $X$,
 so the equaliser is a formal scheme over $X$.

 More explicitly, we have $\XCh(G)=\spf(C(E,G))$, where $C(E,G)$ is
 defined as follows.  We start with $\OX$ and adjoin power series
 variables $c_{ik}$ for $i=1,\ldots,h$ and $k=1,\ldots,d_i$, and put
 $c_{i0}=1$.  We then put $f_i(t)=\sum_{k=0}^{d_i}c_{ik}t^{d_i-k}$ and
 impose the relations obtained by equating coefficients in the
 following identities between polynomials:
 \begin{align*}
  \sum_{a=0}^{d_id_j} p_{d_i,d_j,a}(c_{i*},c_{j*}) t^{d_id_j-a} &= 
   \prod_k f_k(t)^{m_{ijk}} \\
  \sum_{a=0}^{\bcf{d_i}{r}} q_{d_i,r,a}(c_{i*}) t^{\bcf{d_i}{r}-a} &=
   \prod_j f_j(t)^{l^r_{ij}} 
 \end{align*}
 The resulting quotient ring is $C(E,G)$.  
\end{proof}

We next explain how to compare $\XCh(G)$ to $X(G)$.  Let $\CGb$ be the
category whose objects are Lie groups, and whose morphisms are the
conjugacy classes of continuous homomorphisms.  We then have a natural
map
\[ R^+(G) = \coprod_d\CGb(G,U(d)) \xra{B}
            \hbT(BG,\coprod_dBU(d)) \xra{\spf(E^0(-))}
            \hX_X(X(G),\Div^+(\GG)).
\]
By taking adjoints, we obtain a map
$X(G)\xra{}\Map(R^+(G),\Div^+(\GG))$, and one checks easily that this
actually lands in the subscheme
$\XCh(G)\subset\Map(R^+(G),\Div^+(\GG))$ of $\Lm$-semiring
homomorphisms.  We thus have a natural map
\[ \tht_G \: X(G) \xra{} \XCh(G). \]
In terms of our explicit description of $C(E,G)$, the map
$\tht^*\:C(E,G)\xra{}E^0BG$ sends $c_{ik}$ to the $k$'th Chern class
of the representation $V_i$.

It is natural to ask whether a homomorphism
$f\:R(G)\xra{}\Div(\GG)(A)$ of $\Lm$-rings is automatically positive.
We next show that we always have
$f(R_1^+(G))\sse\Div_1^+(\GG)\simeq\GG$, but the corresponding claim
for $d>1$ seems to be false.

\begin{proposition}\label{prop-Div-one}
 If $D\in\Div_1(\GG)(A)$ and $\lm^k(D)=0$ for all $k>1$ then
 $D\in\Div_1^+(\GG)(A)$.
\end{proposition}
\begin{proof}
 We can write $D=E-e[0]$ for some $e\geq 0$ and
 $E\in\Div_{e+1}^+(\GG)(A)$.  Put $D'=\lm^{e+1}E\in\Div_1^+(\GG)$.  We
 have $E=D+e[0]$ so 
 \[ D'=\sum_{i+j=e+1}\lm^i(D)\lm^j(e[0])=\lm^1(D)\lm^e(e[0])=D, \]
 so $D\in\Div_1^+(\GG)$ as claimed.
\end{proof}
\begin{corollary}\label{cor-Div-one}
 If $L\in R_1^+(G)$ and $f\:R(G)\xra{}\Div(\GG)(A)$ is a map of
 $\Lm$-rings then $f(L)\in\Div_1^+(\GG)$.
\end{corollary}
\begin{proof}
 Clearly $\lm^kL=0$ for $k>1$ so $\lm^kf(L)=0$ for $k>1$, and
 $f(L)\in\Div_1(\GG)(A)$ by Lemma~\ref{lem-psipv} so
 $f(L)\in\Div_1^+(\GG)$ by the proposition.
\end{proof}

\begin{proposition}
 For suitable formal groups $\GG$ and rings $A$, there exist divisors
 $D\in\Div_2(\GG)(A)$ such that $\lm^kD=0$ for $k>2$ but
 $D\not\in\Div_2^+(\GG)$.
\end{proposition}
\begin{proof}
 We will assume that $\OX=\F_2$, so $p=2$.  Suppose that
 $a,b\in\GG(A)$ and $2a=2b=0$.  Put $c=a+b$ so $2a=2b=2c=a+b+c=0$, and
 put $E=[a]+[b]+[c]$ and $D=E-[0]$.  Then
 $\lm^2E=[a+b]+[b+c]+[c+a]=[c]+[a]+[b]=E$ and $\lm^3E=[0]=1$ so
 $\lm_t(E)=1+tE+t^2E+t^3=(1+tD+t^2)(1+t)$, so $\lm_t(D)=1+tD+t^2$.
 Thus $\lm^kD=0$ for $k>2$.  If $D$ is in $\Div_2^+(\GG)(A)$ we must
 have $x(a)x(b)x(c)=c_3(E)=c_3(D+[0])=0$.  Note also that
 $x(c)=x(a-b)=x(a)-_Fx(b)$, which is a unit multiple of $x(a)-x(b)$,
 so the condition is equivalent to $x(a)^2x(b)=x(a)x(b)^2$.  The
 universal example for $A$ is
 $\OX\psb{y,z}/([2](y),[2](z))=\F_2[y,z]/(y^{2^n},z^{2^n})$ (where
 $y=x(a),z=x(b)$).  Clearly in this case we have $y^2z\neq yz^2$ so
 $D\not\in\Div_2^+(\GG)$.
\end{proof}

\section{The group $\Sg_3$}\label{sec-Sigma-three}

In this section we work through the case where $G=\Sg_3$ and $E$ is
the $2$-periodic version of Morava $K$-theory at the prime $3$ with
height $2$.  Many constructions discussed here will be generalised
later.  Recall that the coefficient ring is $E^*=\F_3[u,u^{-1}]$,
where $|u|=-2$.

We have a coordinate $x$ on $\GG$ such that 
\begin{align*}
 x(-a)  &= [-1](x(a)) = -x(a) \\ 
 x(3a)  &= [3](x(a)) = x(a)^9 \\
 x(a+b) &= x(a) +_F x(b) = x(a) + x(b) \pmod{x(a)^3x(b)^3}
\end{align*}
for all $a,b\in\GG$.  (The first equation is true because the formal
group law $F$ associated to $E$ has an integral lift whose logarithm
$\log_F(x)=\sum_kx^{9^k}/3^k$ satisfies $\log_F(-x)=-\log_F(x)$.  The
second is well-known, and the third follows from~\cite[Lemma
80]{st:fsf}.) 

Define $y,z\:\Div_2^+(\GG)\xra{}\aff^1$ by $y([a]+[b])=x(a)x(b)$ and
\[ z([a]+[b])=x(\lm^2([a]+[b]))=x(a+b)=x(a)+_Fx(b) \] 
(which is a unit multiple of $x(a)+x(b)$).  One checks that that
$\O_{\Div_2^+(\GG)}=\F_3\psb{y,z}$.  If we let $Z=\SDiv_2^+(\GG)$ be
the scheme of divisors $D\in\Div_2^+(\GG)$ such that $\lm^2(D)=[0]$
then it follows that $\O_Z=\F_3\psb{y,z}/z=\F_3\psb{y}$.  There is an
evident map $\dl\:\GG\xra{}Z$ defined by $\dl(b)=[b]+[-b]$, and
$y(\dl(b))=x(b)x(-b)=-x(b)^2$ so the map
$\dl^*\:\F_3\psb{y}\xra{}\F_3\psb{x}$ sends $y$ to $-x^2$.  In
particular, we see that $\dl$ is finite and faithfully flat, with
degree two.

Next, note that 
\[ \dl(b)^2 = [2b] + [-2b] + 2[0] = \psi^2(\dl(b)) + 2[0]; \]
as $\dl$ is faithfully flat, it follows that $D^2=\psi^2D+2[0]$ for
any $D\in Z$.

Let $Y$ be the scheme of divisors $D\in Z$ such that $\psi^2(D)=D$.
To analyse this, note that
\[ x(2b) = x(-b+3b) = [-1](x(b))+_F[3](x(b)) = 
   -x(b) + x(b)^9 \pmod{x(b)^{12}}, 
\]
so 
\[ -x(2b)^2 = -x(b)^2 - x(b)^{10} \pmod{x(b)^{12}}, \]
or in other words
\[ y(\psi^2\dl(b)) = y(\dl(b)) - y(\dl(b))^5 \pmod{y(\dl(b))^6}. \]
As $\dl$ is faithfully flat, we deduce that
\[ y(\psi^2(D)) = y(D) - y(D)^5 \pmod{y(D)^6} \]
for all $D\in Z$.  It follows that $(\psi^2)^*y-y$ is a unit
multiple of $y^5$ in $\F_3\psb{y}$ and thus that $\O_Y=\F_3[y]/y^5$. 

The character table of $G=\Sg_3$ is
\[ \begin{array}{|c|c|c|c|}
 \hline
      &  1 & \ep & \sg \\ \hline
  1^3 &  1 &   1 &   2 \\ \hline
  1.2 &  1 &  -1 &   0 \\ \hline
  3   &  1 &   1 &  -1 \\ \hline
\end{array}\]
From this we see that 
\[ R(G)=\Z[\ep,\sg]/(\ep^2-1,\ep\sg-\sg,\sg^2-\sg-1-\ep). \]
The only interesting $\lm$-operation is that $\lm^2(\sg)=\ep$.

Let $f\:R^+(G)\xra{}\Div^+(\GG)(A)$ be a $\Lm$-semiring homomorphism,
in other words a point of $\XCh(G)$.  As $\Div_1^+(\GG)\simeq\GG$,
there is a unique point $a\in\GG$ such that $f(\ep)=[a]$.  We also
write $D=f(\sg)\in\Div_2^+(\GG)$.  As $f$ is a map of $\Lm$-semirings,
these satisfy
\begin{align*}
 [2a]   &= [a]^2 = f(\ep^2) = f(1) = [0] \\
 [a]D   &= f(\ep\sg) = f(\sg) = D \\
 D^2    &= f(\sg^2) = f(\sg+1+\ep) = D + [0] + [a] \\
 \lm^2D &= f(\lm^2\sg) = f(\ep) = [a].
\end{align*}
As we work mod $3$, the map $2\:\GG\xra{}\GG$ is an isomorphism so the
first equation gives $a=0$, so the second equation is automatic and
the last equation says that $D\in Z$.  Thus, the third equation
becomes 
\[ \psi^2D + 2[0] = D^2 = D + 2[0]. \]
The semiring $\Div^+(\GG)$ embeds in the ring $\Div(\GG)$ so we can
cancel to see that $\psi^2(D)=D$, so $D\in Y$.  We can thus define a
map $\chi\:\XCh(G)\xra{}Y$ by $\chi(f)=f(\sg)$.  One can check that
the whole argument is reversible, so $\chi$ is an isomorphism and 
\[ C(E,G)=\O_{\XCh(G)}\simeq \O_Y = \F_3[y]/y^5. \]

We also have a short exact sequence $C_3\xra{}G\xra{}C_2$ leading to
an Atiyah-Hirzebruch-Serre spectral sequence
\[ H^p(C_2;E^qBC_3) \convto E^{p+q}BG. \]
We have $E^*BC_3=\F_3[u^{\pm 1}][x]/x^9$ (where $|u|=2$ and $|x|=0$),
and $C_2$ acts on this by $u\mapsto u$ and $x\mapsto [-1](x)=-x$.
Because $C_2$ has order coprime to $3$ we see that the spectral
sequence collapses to an isomorphism
\[ E^*BG = (E^*BC_3)^{C_2} = E^*[y]/y^5, \]
where $y=-x^2$.  After some comparison of definitions we see that the
map $\tht_G\:X(G)\xra{}\XCh(G)$ is an isomorphism.

\section{The Abelian case}\label{sec-Abelian}

\begin{theorem}\label{thm-Abelian}
 If $G$ is Abelian then $\tht_G\:X(G)\xra{}\XCh(G)$ is an
 isomorphism. 
\end{theorem}
\begin{proof}
 Put $G^*=\Hom(G,S^1)$, so $R^+(G)=\N[G^*]$, and let $A$ be an
 $\OX$-algebra.  If $f\:\N[G^*]\xra{}\Div^+(\GG)(A)$ is a
 $\Lm$-semiring homomorphism, then $f$ induces a group homomorphism
 $f'\:G^*=R^+_1(G)\xra{}\GG(A)=\Div^+_1(\GG)(A)$.  Conversely, given a
 group homomorphism $f'\:G^*\xra{}\GG(A)$ we get a map
 $R^+(G)=\N[G^*]\xra{}\N[\GG(A)]$ of $\Lm$-semirings.  We can compose
 this with the map $\N[\GG(A)]\xra{}\Div^+(\GG)(A)$ in the proof of
 Proposition~\ref{prop-Div} to get a map $R^+(G)\xra{}\Div^+(\GG)(A)$,
 or in other words a point of $\XCh(G)(A)$.  One checks that these
 constructions give a bijection $\Hom(G^*,\GG(A))\simeq\XCh(G)(A)$, or
 equivalently an isomorphism $\Hom(G^*,\GG)\simeq\XCh(G)$.  There is
 also an isomorphism $X(G)\simeq\Hom(G^*,\GG)$ (see~\cite[Proposition
 2.9]{grst:vlc}), so we have an isomorphism $X(G)\simeq\XCh(G)$.  A
 straightforward comparison of definitions shows that this is the same
 as $\tht_G$.
\end{proof}

\section{The height one case}\label{sec-height-one}

In this section we choose a prime $p$ and let $E$ be the $p$-complete
complex $K$-theory spectrum.  We thus have $X=\spf(\Zp)$, so 
discrete $\OX$-algebras are just $p$-torsion rings. We also have
$\GG=\Gmh$, so
\[ \GG(A) = \{u\in A^\tm \st 1-u \text{ is nilpotent } \}. \]

\begin{theorem}
 If $E$ is the $p$-adic complex $K$-theory spectrum and $G$ is a
 $p$-group then the map $\tht_G\:X(G)\xra{}\XCh(G)$ is an isomorphism.
\end{theorem}

The rest of this section constitutes the proof.  We fix a $p$-group
$G$ and write $\tht=\tht_G$ for brevity.

The first ingredient is the Atiyah-Segal completion theorem.  In the
case of a $p$-group, this says that
\[ \O_{X(G)} = E^0BG = \Zp\ot R(G). \]

We know that $\Div(\Gmh)$ is the free ring scheme generated by the
group scheme $\Gmh$.  Recall that the affine line $\aff^1$ is just the
forgetful functor from $p$-torsion rings to sets.  This is a ring
scheme in a natural way, and it contains $\Gmh$ as a subgroup of its
group of units.  We thus have a ring map
\[ \xi\:\Div(\Gmh)\xra{}\aff^1 \] 
extending the inclusion of $\Gmh$ in $\aff^1$.  If
$D=\sum_in_i[u_i]\in\Div(\Gmh)(A)$ then
$\xi(D)=\sum_in_iu_i\in\aff^1(A)=A$. 

Next recall that $\Div_d^+(\Gmh)=\Gmh^d/\Sg_d$, so to describe a
function $f\:\Div_d^+(\Gmh)\xra{}\aff^1$ it suffices to give the
symmetric function $\ov{f}\:\Gmh^d\xra{}\aff^1$ such that
\[ f(\sum_{i=1}^d [u_i]) = \ov{f}(u_1,\ldots,u_d). \]
Let $\sg_j\:\aff^d\xra{}\aff^1$ be the $j$'th elementary symmetric
function and define 
\begin{align*}
  c_j(\sum_i[u_i]) &= \sg_j(1-u_1,\ldots,1-u_d)         \\
  a_j(\sum_i[u_i]) &= \sg_j(u_1,\ldots,u_d)             \\
 a'_j(\sum_i[u_i]) &= a_j(\sum_i[u_i]) - \bcf{d}{j}.
\end{align*}
Recall that $\O_{\Div_d^+(\Gmh)}$ is the set of all maps
$\Div_d^+(\Gmh)\xra{}\aff^1$, so $c_j$, $a_j$ and $a'_j$ can be viewed
as elements of this ring.  The function $u\mapsto 1-u$ is a coordinate
on the formal group $\Gmh$ so a well-known argument gives an
isomorphism 
\[ \O_{\Div_d^+(\Gmh)} = \Zp\psb{c_1,\ldots,c_d}. \]
It is not hard to deduce that
\[ \O_{\Div_d^+(\Gmh)} = \Zp\psb{a'_1,\ldots,a'_d}, \]
which is the completion of the ring $\Zp[a_1,\ldots,a_d]$ at the ideal
generated by the elements $a'_i=a_i-\bcf{d}{i}$.  Note also that
$a_1(D)=\xi(D)$ for $D\in\Div_d^+(\Gmh)$.

\begin{lemma}
 For any divisor $D\in\Div_d^+(\Gmh)$ we have $a_j(D)=\xi(\lm^j(D))$.
\end{lemma}
\begin{proof}
 We may assume that $D=\sum_i[u_i]$ for some elements
 $u_1,\ldots,u_d\in\Gmh$.  For any $I\sse\{1,\ldots,d\}$ with $|I|=j$
 we put $u_I=\prod_{i\in I}u_i$.  We then have $\lm^jD=\sum_I[u_I]$
 and thus 
 \[ \xi\lm^j(D)=\sum_Iu_I=\sg_j(u_1,\ldots,u_d)=a_j(D). \]
\end{proof}

Let $V$ be a complex vector bundle over a space $Z$ with associated
projective bundle $PV$, and let $\DD(V)=\spf(E^0PV)$ be the
corresponding divisor on $\GG$ over $\spf(E^0Z)$.  If
$V=\bigoplus_{i=1}^dL_i$ for some complex line bundles $L_i$ then each
$L_i$ can be regarded as an element of $E^0Z=K^0(Z;\Zp)$, and one sees
easily that $\DD(V)=\sum_i[L_i]$ so
$a_j(\DD(V))=\sg_j(L_1,\ldots,L_d)=\lm^j(V)$.  By the splitting
principle we see that 
\[ \xi\lm^j(\DD(V)) = a_j(\DD(V)) = \lm^j(V) \]
even when $V$ does not split as a sum of line bundles.  

Now consider the case $Z=BG$ and suppose that $V$ comes from a
representation of $G$, which we also call $V$.  Suppose we have a
point $x\in X(G)(A)$ for some $p$-torsion ring $A$, corresponding to a
ring map $\hx\:E^0BG\xra{}A$.  One can now check from the definitions
that $\tht(x)(V)=\hx_*(\DD(V))$, so $f(\tht(x)(V))=\hx(f(\DD(V)))$ for
any $f\in\O_{\Div_d^+(\Gmh)}$.  In particular, we have
\[ a_j(\tht(x)(V)) = \hx(a_j(\DD(V))) = \hx(\lm^j(V)). \]

We can now construct the map $\zt\:\XCh(G)\xra{}X(G)$ that will turn
out to be inverse to $\tht$.  Let $A$ be a $p$-torsion ring.  A positive
homomorphism $f\in\XCh(G)(A)$ gives rise to a homomorphism
$\xi_A\circ f\:R(G)\xra{}\aff^1(A)=A$, which factors canonically
through $\Zp\ot R(G)=E^0BG=\O_{X(G)}$ because $A$ is a $p$-torsion
ring.  This gives a continuous homomorphism $\O_{X(G)}\xra{}A$, or in
other words a point of $X(G)(A)$, which we call $\zt(f)$.  This
construction gives a natural map $\zt\:\XCh(G)\xra{}X(G)$, as
required.

Suppose we start with a point $x\in X(G)(A)$, corresponding to a ring
map $\hx\:\Zp\ot R(G)=E^0BG\xra{}A$.  We need to check that
$\zt\tht(x)=x$, or equivalently that $\xi_A(\tht(x)(V))=\hx(V)$ for
all $V\in \Zp\ot R(G)$, and it will suffice to do this for all honest
representations $V\in R_d^+(G)$ for all $d$.  In that context we have
$\xi=a_1$ so
\[ \xi_A(\tht(x)(V)) = a_1(\tht(x)(V)) = \hx(\lm^1(V)) = \hx(V)
\]
as required.  Thus $\zt\tht=1_{X(G)}$.

Suppose instead that we start with a point $f\in\XCh(G)(A)$, in other
words a positive homomorphism $f\:R^+(G)\xra{}\Div^+(\Gmh)(A)$.  We
need to check that $(\tht(\zt(f)))(V)=f(V)\in\Div_d^+(\Gmh)(A)$ for
all $V\in R^+_d(G)$.  To see this, note that
\begin{align*}
 a_j(\tht(\zt(f))(V)) &= \widehat{\zt(f)}(\lm^j(V)) \\
  &= \xi(f(\lm^j(V))) \\
  &= \xi(\lm^j(f(V))) \\
  &= a_j(f(V)).
\end{align*}
It follows that $a'_j(\tht(\zt(f))(V))=a'_j(f(V))$ and the functions
$a'_j$ generate $\O_{\Div_d^+(\Gmh)}$ so $\tht(\zt(f))(V)=f(V)$, as
required.  This shows that $\tht\zt=1_{\XCh(G)}$, so $\tht$ is an
isomorphism as claimed.

\section{Reduction to the Sylow subgroup}\label{sec-Sylow}

By a well-known transfer argument, if $P$ is a Sylow $p$-subgroup of
$G$ then the restriction map $E^0BG\xra{}E^0BP$ is injective, and
similar methods give some control over the image.  In this section we
develop some analogous results for the approximation $C(E,G)$ to
$E^0BG$.  Let $I$ be the kernel of the restriction map
$R(G)\xra{}R(P)$ (which is independent of the choice of $P$).  Note
that $R(G)/I$ is isomorphic to the image of the restriction map, so it
is a free Abelian group of finite rank and it inherits a $\Lm$-ring
structure.

\begin{proposition}\label{prop-RGI}
 Any $\Lm$-ring homomorphism $f\:R(G)\xra{}\Div(\GG)(A)$ factors
 through $R(G)/I$.
\end{proposition}
\begin{proof}
 As usual, we let the exponent of $G$ be $e=p^ve'$, where $e'$ is
 coprime to $p$.  If $V\in I$ and $g\in G$ has $p$-power order then
 $g$ is conjugate to an element of $P$ and thus $\chi_V(g)=0$.  If $g$
 is an arbitrary element of $g$ then the order of $g$ divides $p^ve'$
 so the order of $g^{e'}$ divides $p^v$, so $\chi_V(g^{e'})=0$.  This
 proves that $\psi^{e'}(V)=0$, so $\psi^{e'}(f(V))=0$ in
 $\Div(\GG)(A)$.  However, the action of $\psi^{e'}$ on $\Div(\GG)(A)$
 is induced by the action of $e'$ on $\GG$, which is invertible because
 $e'$ is coprime to $p$.  This implies that $f(V)=0$ as claimed.
\end{proof}

\begin{corollary}
 If $|G|$ is prime to $p$ then $\XCh(G)=X(G)=X$. \qed
\end{corollary}

\begin{corollary}\label{cor-RGI}
 $\XCh(G)$ is a closed subscheme of the scheme of $\Lm$-ring maps
 $R(G)/I\xra{}\Div(G)$. \qed
\end{corollary}

\begin{remark}\label{rem-restriction}
 The disadvantage of this point of view is that the positivity
 conditions $f(R^+(G))\sse\Div^+(\GG)$ become less visible when we
 work with $R(G)/I$.  However, the situation simplifies again if we
 assume that $P$ is normal in $G$.  In that case, it is known that the
 map $R^+(G)\xra{}R^+(P)^G$ is surjective; this was first proved by
 Gallagher~\cite{ga:gcn}, and we will give an alternative proof in
 Section~\ref{sec-restriction}.  The sums of $G$-orbits of
 irreducibles in $R^+(P)$ give a canonical system of generators for
 $R^+(P)^G$, which we can lift to get representations
 $\rho_1,\ldots,\rho_t$ of $G$ say.  Let $d_i$ be the degree of
 $\rho_i$.  If $\rho\in R_d^+(G)$ then $\res^G_P(\rho)\in R_d^+(P)^G$
 so $\res^G_P(\rho)\simeq\sum_im_i\res^G_P(\rho_i)$ for some integers
 $m_i\geq 0$ with $\sum_im_id_i=d$, so $\rho-\sum_im_i\rho_i\in I$.
 Thus, for any map $f\:R(G)\xra{}\Div(\GG)$ of $\Lm$-rings we have
 $f(\rho)=\sum_im_if(\rho_i)$, so if $f(\rho_i)\in\Div_{d_i}^+(\GG)$
 for all $i$ then $f(\rho)\in\Div_d^+(\GG)$.  Using this, we see that
 $\XCh(G)$ is the scheme of $\Lm$-ring maps $R(G)/I\xra{}\Div(\GG)$
 such that $f(\rho_i)\in\Div_{d_i}^+(\GG)$ for $i=1,\ldots,t$.
\end{remark}

We next give two results that help us to understand $R(G)/I$ without
computing $R(G)$.

\begin{proposition}
 Let $h$ be the number of conjugacy classes of elements $g\in G$ whose
 order is a power of $p$.  Then $R(G)/I\simeq\Z^h$ as Abelian groups.
\end{proposition}
\begin{proof}
 We already know that $R(G)/I$ is a free Abelian group, so we just
 need to determine its rank, so it is enough to show that
 $\C\ot R(G)/I\simeq\C^h$.  Let $C$ be the set of conjugacy classes of
 $p$-power order, and let $C'$ be the set of all other conjugacy
 classes.  Let $U(G)$ be the ideal of virtual representations
 $V\in R(G)$ whose character is zero on $C'$.  It is well-known that
 $\C\ot R(G)\simeq F(C\amalg C',\C)$.  This isomorphism carries $U(G)$
 to $F(C,\C)$ and $I$ to $F(C',\C)$ so the map
 $\C\ot U(G)\xra{}\C\ot R(G)/I$ is an isomorphism.  Clearly
 $\dim_\C\C\ot U(G)=|C|=h$, and the claim follows.
\end{proof}
\begin{remark}
 The proof shows that $U(G)$ is a subgroup of finite index in
 $R(G)/I$.  This index need neither be a power of $p$ nor coprime to
 $p$, as one sees by taking $G=\Sg_3$ and $p=2$ or $3$; the index is
 $2$ in both cases.
\end{remark}

\begin{proposition}
 There is a natural isomorphism $\Zp\ot R(G)/I=KU^0_p(BG)$ (where
 $KU_p$ is the $p$-adic completion of the complex $K$-theory
 spectrum.) 
\end{proposition}
\begin{proof}
 Let $KU_G$ be the usual $G$-spectrum for equivariant $K$-theory, and
 let $L_G$ be its $p$-completion.  It is well-known that
 $KU^0_G(S^0)=R(G)$, which is free Abelian of finite rank, and it
 follows that $L_G^0(S^0)=R(G)^\wedge_p=R(G)\ot\Zp$.  We give this
 ring and all its quotients the $p$-adic topology, or equivalently the
 profinite topology, which is compact.  The argument of the
 Atiyah-Segal completion theorem shows that $KU_p^0BG=L^0EG$ is the
 completion of $R(G)^\wedge_p$ at the augmentation ideal $J_G$.  By a
 compactness argument, we deduce that the map
 $R(G)^\wedge_p\xra{}KU_p^0BG$ is surjective; the kernel is
 $J_G^\infty:=\bigcap_kJ_G^k$.  Now let $P$ be a Sylow $p$-subgroup,
 so by the same arguments $KU_p^0BP=R(P)^\wedge_p/J_P^\infty$.  It is
 well-known that $J_P^N\leq pR(P)^\wedge_p$ for $N\geq 0$ (use the
 fact that $x^p-\psi^p(x)\in pR(P)$ for all $x\in R(P)$, for example)
 and it follows that $J_P^\infty=0$, so $KU^p_0BP=R(P)^\wedge_p$.  We
 now have a diagram as follows.
 \begin{diag}
  \node[2]{J_G^\infty} \arrow{e}   \arrow{s,V}
  \node{J_P^\infty=0}              \arrow{s}           \\
  \node{I^\wedge_p}    \arrow{e,V} 
  \node{R(G)^\wedge_p} \arrow{e}   \arrow{s,A}
  \node{R(P)^\wedge_p}             \arrow{s,r}{\simeq} \\
  \node[2]{KU_p^0BG}   \arrow{e,V} \node{KU^0_pBP.}
 \end{diag}
 We have seen that the columns are short exact.  As $\Zp$ is flat over
 $\Z$ and $I$, $R(G)$ and $R(P)$ are finitely generated Abelian
 groups, we see that the middle row is left exact.  The bottom
 horizontal map is injective by a transfer argument.  By a diagram
 chase we deduce that $I^\wedge_p=J_G^\infty$, so
 $KU^0_pBG=R(G)^\wedge_p/I^\wedge_p=(R(G)/I)\ot\Zp$ as claimed.
\end{proof}

\section{Finiteness}\label{sec-finite}

It is a fundamental fact that the scheme $X(G)$ is finite over $X$, or
equivalently that $E^0BG$ is a finitely generated module over $E^0$.  
This is proved in the present generality 
as~\cite[Corollary 4.4]{grst:vlc}; the argument is the same as
in~\cite{hokura:ggc}.  In this section we show that $\XCh(G)$ is also
finite over $X$.  We also study some auxiliary schemes that come up in
the proof, as they turn out to be useful in specific computations.

\begin{definition}
 For any $v\geq 0$ we put
 \[ \GG(v) = \ker(p^v\:\GG\xra{}\GG). \]
 In terms of our coordinate, we have
 $\GG(v)=\spf(E^0\psb{x}/[p^v](x))$.  Next, recall that
 $\GG(v)\leq\GG$ is a divisor of degree $p^{nv}$, where $n$ is the
 height of $\GG$.  For any $m\geq 0$ we let $\GG(v,m)$ be $m$ times
 this divisor, considered as a subscheme of $\GG$; in other words
 $\GG(v,m)=\spf(E^0\psb{x}/[p^v](x)^m)$.  For any formal scheme $Y$
 over $X$ such that $\O_Y$ is a free module of finite rank $r$ over
 $\OX$, we define $Y^d/\Sg_d$ to be $\spf$ of the $d$'th symmetric
 tensor power of $\O_Y$ over $\O_X$.  If $\{e_1,\ldots,e_r\}$ is a
 basis for $\OY$ over $\OX$ then the monomials
 $e^\al:=\prod_{i=1}^re_i^{\al_i}$ with $\sum_i\al_i=d$ form a basis
 for $\O_{Y^d/\Sg_d}$, so this ring is free of rank $\bcf{r+d-1}{d}$
 over $\O_X$.  We will use this construction in the cases $Y=\GG(v)$
 and $Y=\GG(v,m)$.  Finally, we define $Z(v,d)$ to be the scheme of
 divisors $D\in\Div_d^+(\GG)$ such that $\psi^{p^v}D=d[0]$.
\end{definition}

\begin{theorem}\label{thm-Zdv}
 The scheme $Z(d,v)$ is finite and flat over $X$, of degree $p^{ndv}$.
 There are closed inclusions
 \[ \GG(v)^d/\Sg_d \xra{} Z(v,d) \xra{} \GG(v,d)^d/\Sg_d
      \xra{} \Div_d^+(\GG).
 \]
 The first two of these are infinitesimal thickenings, in other words
 the corresponding maps of rings are surjective with nilpotent
 kernel.  If $\OX$ is a field (necessarily of characteristic $p$) then
 $Z(v,d)$ is the fibre of the $nv$-fold relative Frobenius map
 \[ F^{nv}_{\Div_d^+(\GG)/X}\:
     \Div_d^+(\GG)\xra{}\Div_d^+((F_X^{nv})^*\GG),
 \]
 so 
 \[ \O_{Z(v,d)} = \OX\psb{c_1,\ldots,c_d}/(c_i^{p^{nv}}). \]
\end{theorem}
\begin{proof}
 First suppose that $\OX$ is a complete regular local ring.  (In the
 topological context, this occurs when $E$ is Landweber exact.)
 Consider the following diagram:
 \begin{diag}
  \node{Z(v,d)} \arrow{s} \arrow{e,V}
  \node{\Div_d^+(\GG)} \arrow{s,l,A}{\psi^{p^v}}
  \node{\GG^d} \arrow{s,r,A}{p^v} \arrow{w,t,A}{\pi} \\
  \node{X} \arrow{e,b}{\zt} 
  \node{\Div_d^+(\GG)}
  \node{\GG^d} \arrow{w,b,A}{\pi}
 \end{diag}
 In the right hand square, all the corresponding rings are complete
 regular local rings.  A finite injective map of such rings always
 makes the target into a free module over the source~\cite[2.2.7 and
 2.2.11]{brhe:cmr}.  The maps $\pi^*$ and $(p^v)^*$ are finite
 injective maps of degrees $d!$ and $p^{ndv}$.  It follows that
 $(\psi^{p^v})^*$ is finite and injective, and thus (as
 $\deg(fg)=\deg(f)\deg(g)$ in this context) that $\psi^{p^v}$ is flat
 of degree $p^{ndv}$.  The left hand square is a pullback by
 definition, and it follows that $Z(v,d)$ is flat of degree $p^{ndv}$
 over $X$.  Using~\cite[Proposition 5.2]{grst:vlc}, it is not hard to
 deduce that this result remains true even if $\OX$ is not regular.

 Next, let $x$ be a coordinate on $\GG$.  Then $\{x^i\st i<p^{nv}\}$
 is a basis for $\O_{\GG(v)}$ over $\OX$, and $\{x^i\st i<p^{nmv}\}$
 is a basis for $\O_{\GG(v,m)}$.  Using this we obtain bases for the
 rings 
 \[ A = \O_{\GG(v)^d} = \OX\psb{x_1,\ldots,x_d}/([p^v](x_i)) \] 
 and
 \[ A' = \O_{\GG(v,m)^d} = \OX\psb{x_1,\ldots,x_d}/([p^v](x_i)^m) \] 
 that are permuted by $\Sg_d$, and the orbit sums give bases for the
 rings 
 \[ B = \O_{\GG(v)^d/\Sg_d} = A^{\Sg_d} \]
 and
 \[ B' = \O_{\GG(v,m)/\Sg_d} = (A')^{\Sg_d}. \]
 Using these, it is easy to see that the map $B'\xra{}B$ is
 surjective, so the map $\GG(v)^d/\Sg_d\xra{}\GG(v,m)^d/\Sg_d$ is a
 closed inclusion.  A similar argument shows that $\GG(v,m)^d/\Sg_d$
 is a closed subscheme of $\Div_d^+(\GG)$.

 Next, put 
 \[ J = \ker(A'\xra{}A) = ([p^v](x_1),\ldots,[p^v](x_d)). \]
 This is clearly a nilpotent ideal, and $\ker(B'\xra{}B)=B'\cap J$
 which is a nilpotent ideal in $B'$.  Thus our map
 $\GG(v)^d/\Sg_d\xra{}\GG(v,m)^d/\Sg_d$ is an infinitesimal
 thickening. 

 It is clear that $\GG(v)^d/\Sg_d$ is contained in $Z(v,d)$.  Next,
 let $W(v,d)$ be the preimage of $Z(v,d)$ in $\GG^d$, or equivalently
 the scheme of $d$-tuples $\ua=(a_1,\ldots,a_d)\in\GG^d$ such that
 $\sum_i[p^va_i]=d[0]$.  If $\ua\in W(v,d)$ then for each $i$ we have
 $p^va_i\in d[0]$ so $a_i\in d.(p^v)^{-1}[0]=d.\GG(v)=\GG(v,d)$.  This
 means that $\ua\in\GG(v,d)^d$ and thus
 $\sum_i[a_i]\in\GG(v,d)^d/\Sg_d$, so the map
 $W(v,d)\xra{\pi}Z(d,v)\xra{}\Div_d(\GG)$ factors through
 $\GG(v,d)^d/\Sg_d$.  As $\pi$ is faithfully flat, it follows that
 $Z(v,d)\sse\GG(v,d)^d/\Sg_d$ as claimed.

 We now have maps
 \[ \GG(v)^d/\Sg_d \xra{i} Z(v,d) \xra{j} \GG(v,d)^d/\Sg_d 
      \xra{k} \Div_d^+(\GG).
 \]
 We know that $ji$, $k$ and $kj$ are closed inclusions and that $ji$
 is an infinitesimal thickening.  It follows easily that $i$, $j$ and
 $k$ are closed inclusions and $i$ and $j$ are infinitesimal
 thickenings. 

 Now suppose that $X$ is a field of characteristic $p$.  We then have
 an iterated Frobenius map $F^{nv}_X\:X\xra{}X$ corresponding to the
 ring map $a\mapsto a^{p^{nv}}$ and thus a formal group
 $\GG'=(F^{nv}_X)^*\GG$ over $X$.  The map $F^{nv}_\GG$ gives rise to
 a map $f=F^{nv}_{\GG/X}\:\GG\xra{}\GG'$.  As $\GG$ has height $n$,
 the map $p^v\:\GG\xra{}\GG$ factors as $\GG\xra{f}\GG'\xra{g}\GG$,
 where $g$ is an isomorphism.  This is just the geometric statement of
 the fact that $[p^v](x)=\gm(x^{p^{nv}})$ for some invertible power
 series $\gm$.  By definition $Z(v,d)$ is the fibre of the map
 $\Div_d^+(\GG)\xra{}\Div_d^+(\GG)$ induced by $p^v\:\GG\xra{}\GG$,
 and it follows easily that it is also the fibre of the map
 $\Div_d^+(\GG)\xra{}\Div_d^+(\GG')$ induced by $f$.  It is easy to
 identify this with the map $F^{nv}_{\Div_d^+(\GG)}$.  If we use the
 usual generators for the coordinate rings of $\Div_d^+(\GG)$ and
 $\Div_d^+(\GG')$ then the corresponding ring map sends $c_k$ to
 $c_k^{p^{nv}}$, so
 $\O_{Z(v,d)}=\OX\psb{c_1,\ldots,c_d}/(c_k^{p^{nv}})$. 
\end{proof}

\begin{corollary}
 The scheme $\XCh(G)$ is finite over $X$. 
\end{corollary}
\begin{proof}
 Let $V_1,\ldots,V_h$ be the irreducible representations of $G$, and
 let $d_1,\ldots,d_h$ be their degrees.  As in the proof of
 Proposition~\ref{prop-scheme}, we see that $\XCh(G)$ is a closed
 subscheme of $\prod_i\Div^+_{d_i}(\GG)$.  Now let $p^v$ be the
 $p$-part of the exponent of $G$.  We see from Lemma~\ref{lem-psipv}
 that $\XCh(G)$ is actually contained in $\prod_iZ(v,d_i)$, which is
 finite over $X$ by the theorem.  It follows that $\XCh(G)$ itself is
 finite, as claimed.
\end{proof}

\begin{corollary}\label{cor-psi-period}
 For any divisor $D\in\Div(\GG)(A)$ there exist $w\geq 0$ such that
 $\psi^k(D)=\psi^l(D)$ whenever $k,l\in\Zp$ with $k=l\pmod{p^w}$.
\end{corollary}
\begin{proof}
 We can reduce easily to the case where $D\in\Div_d^+(\GG)$ for some
 $d\geq 0$.  Lemma~\ref{lem-psipv} tells us that $\psi^{p^v}D=d[0]$
 for some $v$, so $D\in Z(v,d)(A)\leq\GG(v,d)^d/\Sg_d(A)$.  Next note
 that $A$ is a discrete $\OX$-algebra so $p$ is nilpotent in $A$, say
 $p^r=0$.  It follows that $[p^r](x)=f(x^p)$ for some power series $f$
 with $f(0)=0$, and thus that $[p^{rs}](x)$ is divisible by $x^{p^s}$.
 Thus, for large $u$ we have $[p^u](x)=0\pmod{x^d}$, so
 $[p^{u+v}](x)=0\pmod{[p^v](x)^d}$, so $\GG(v,d)\leq\GG(u+v)$.  If we
 put $w=u+v$ this tells us that $D\in\GG(w)^d/\Sg_d$, and the action
 of $\psi^k$ on $\GG(w)^d/\Sg_d$ clearly depends only on the
 congruence class of $k$ mod $p^w$, as required. 
\end{proof}

\section{Generalised character theory}\label{sec-gen-char}

In~\cite{hokura:ggc}, Hopkins, Kuhn and Ravenel describe $\Q\ot E^0BG$
in terms of ``generalised characters''.  In this section we will give
an analogous but less precise description of $\Q\ot C(E,G)$.

To explain the HKR theory, write $\Tht=(\QZp)^n$ and
$\Tht^*=\Hom(\Tht,\QZp)=\Z_p^n$.  (Elsewhere these are denoted by
$\Lm$ and $\Lm^*$, but there are enough $\Lm$'s in this paper
already.)  Let $\Tht(v)$ be the subgroup of $\Tht$ killed by $p^v$,
and let $\Level(v,\GG)$ be the scheme of maps $\phi\:\Tht(v)\xra{}\GG$
such that $\sum_{a\in\Tht(v)}[\phi(a)]\leq\GG(v)$ in $\Div^+(\GG)$.
See~\cite{st:fsf} for more information about these schemes.  Put
$D_m=\O_{\Level(m,\GG)}$ and $D=\colim_mD_m$ and $L=\Q\ot D$.  This is
a free module of countable rank over $\Q\ot\OX$.  If $\GG$ is a
universal deformation (as in the case considered by HKR) then it can
be described more explicitly: the Weierstrass preparation theorem
implies that $[p^v](x)$ is a unit multiple of a monic polynomial
$g_v(x)$ of degree $p^{nv}$, and $L$ is obtained from $\Q\ot\OX$ by
adjoining full set of roots for $g_v(x)$ for all $v$.

Now let $\Om(G)$ be the set of $G$-conjugacy classes of homomorphisms
$\Tht^*\xra{}G$, and let $F(\Om(G),L)$ be the ring of all functions
$u\:\Om(G)\xra{}L$ (with pointwise operations).  HKR construct an
isomorphism
\[ \tau \: L\ot_{\OX}E^0BG \xra{} F(\Om(G),L). \]
They work with a particular admissible cohomology theory $E$, but it
is not hard to extend their result to all admissible theories;
see~\cite[Proposition 5.2 and Appendix B]{grst:vlc} for some pointers.

Now consider the $\Lm$-semiring $\N[\Tht]=\coprod_d\Tht^d/\Sg_d$ and
the $\Lm$-ring $\Z[\Tht]$.  If we give $\Tht^*$ its $p$-adic topology
then every subgroup of finite index is open and any continuous
homomorphism $\Tht^*\xra{}GL_n(\C)$ factors through a finite quotient
of $\Tht^*$.  Using this we can identify $\N[\Tht]$ with the semiring
of continuous representations of $\Tht^*$, and $\Z[\Tht]$ with the
corresponding ring of virtual representations.

\begin{definition}
 We say that a $\Lm$-ring homomorphism $f\:R(G)\xra{}\Z[\Tht]$ is
 \emph{positive} if $f(R^+(G))\sse\Z[\Tht]^+$.  We write $\OmCh(G)$
 for the set of $\Lm$-semiring homomorphisms $R^+(G)\xra{}\Z[\Tht]^+$,
 or equivalently the set of positive $\Lm$-ring homomorphisms
 $R(G)\xra{}\Z[\Tht]$.
\end{definition}
\begin{remark}\label{rem-OmCh-i}
 The arguments of Lemma~\ref{lem-psipv}, Corollary~\ref{cor-Div-one}
 and Proposition~\ref{prop-RGI} show that any positive homomorphism
 $f\:R(G)\xra{}\Z[\Tht]$ of $\Lm$-rings automatically sends $R^+_d(G)$
 to $\Z[\Tht]^+_d$ and $I$ to $0$ (where $I$ is the kernel of the
 restriction map to a Sylow subgroup).  There is also an evident
 analogue of Remark~\ref{rem-restriction} for $\OmCh(G)$.
\end{remark}
\begin{remark}\label{rem-OmCh-ii}
 From the definitions we know that the $\Lm$-operations determine the
 Adams operations.  Conversely, it is well-known and easy to check
 that the Adams operations determine the $\Lm$-operations rationally.
 As $\Z[\Tht]$ is torsion-free, it follows that a ring map
 $f\:R(G)\xra{}\Z[\Tht]$ preserves the $\Lm$-operations iff it
 preserves the Adams operations.  The corresponding statement for
 homomorphisms $R(G)\xra{}\Div(\GG)$ is false, however.
\end{remark}

\begin{theorem}\label{thm-gen-char}
 There are natural maps $\kp\:\Om(G)\xra{}\OmCh(G)$ and
 $\tCh\:L\ot C(E,G)\xra{}L\ot E^0BG$ making the following diagram
 commute:
 \begin{diag}
  \node{L\ot C(E,G)} \arrow{s,l}{\tCh} \arrow{e,t}{1\ot\tht^*}
  \node{L\ot E^0BG}  \arrow{s,r}{\tau} \\
  \node{F(\OmCh(G),L)}                 \arrow{e,b}{\kp^*}
  \node{F(\Om(G),L).}
 \end{diag}
 (Here the tensor products are taken over $E^0$.)  Moreover, the map
 $\tCh$ is surjective with nilpotent kernel.
\end{theorem}
\begin{proof}
 For brevity we will write $C(L,G)=L\ot C(E,G)$ and
 $L^0BG=L\ot E^0BG$.  This is a slight abuse because these functors do
 not arise from a spectrum $L$.  We also let $v$ be any integer
 greater than or equal to the $p$-adic valuation of the exponent of
 $G$.

 Any homomorphism $u\:\Tht^*\xra{}G$ factors through
 $\Tht^*/p^v=\Tht(v)^*$ and thus is automatically continuous (for the
 discrete topology on $G$).  It thus gives a positive homomorphism
 $u^*\:R(G)\xra{}\Z[\Tht]$, and it is well-known that this depends
 only on the conjugacy class of $u$, so this construction gives a
 natural map $\kp\:\Om(G)\xra{}\OmCh(G)$.

 It is easy to see using Adams operations that any positive
 homomorphism $f\:R(G)\xra{}\Z[\Tht]$ actually lands in the subring
 $\Z[\Tht(v)]$.  Suppose we have a level structure
 $\phi\:\Tht(v)\xra{}\GG(A)$.  As $\Tht(v)$ is a finite Abelian group,
 this gives rise as in Theorem~\ref{thm-Abelian} to a positive
 homomorphism $R(\Tht^*/p^v)=\Z[\Tht(v)]\xra{}\Div(\GG)(A)$, which we
 can compose with $f$ to get a positive homomorphism
 $R(G)\xra{}\Div(\GG)(A)$, or in other words a point of $\XCh(G)(A)$,
 which we call $\rhoCh(f,\phi)$.  This construction produces a map
 $\rhoCh\:\OmCh(G)\tm\Level(v,\GG)\xra{}\XCh(G)$ of formal schemes
 over $X$, corresponding to a map
 $\rhoCh^*\:C(E,G)\xra{}F(\OmCh(G),D_v)\subset F(\OmCh(G),L)$.  After
 tensoring by $L$ we obtain the required map
 $\tCh\:C(L,G)\xra{}F(\OmCh(G),L)$.

 We next recall the definition of $\tau$.  Suppose that
 $u\:\Tht^*/p^v\xra{}G$ and
 $\phi\in\Level(v,\GG)\subset\Hom(\Tht(v),\GG)=X(\Tht^*/p^v)$.  We
 then have a point $\rho(u,\phi):=X(u)(\phi)\in X(G)$.  This
 construction gives a map $\rho\:\Om(G)\tm\Level(v,\GG)\xra{}X(G)$ and
 thus a map $\rho^*\:E^0BG\xra{}F(\Om(G),D_v)\subset F(\Om(G),L)$.
 After tensoring by $L$ we obtain the required map $\tau$.

 One can check from the definitions that the following diagram
 commutes:
 \begin{diag}
  \node{\XCh(G)} 
  \node{X(G)}
  \arrow{w,t}{\tht} \\
  \node{\OmCh(G)\tm\Level(m,\GG)}
  \arrow{n,l}{\rhoCh}
  \node{\Om(G)\tm\Level(m,\GG).}
  \arrow{w,b}{\kp\tm 1}
  \arrow{n,r}{\rho} 
 \end{diag}
 It follows easily that the diagram in the statement of the theorem
 commutes. 

 To understand $\tCh$ more explicitly, let $\phi\in\Level(v,\GG)(D_v)$
 be the universal example of a level structure.  For any element
 $a\in\Tht(v)$ we then have a point $\phi(a)\in\GG(D_v)$ and thus an
 element $x_a:=x(\phi(a))\in D_v$.  These elements satisfy
 $x_{a+b}=x_a+_Fx_b$ and $\GG(v)=\sum_{a\in\Tht(v)}[\phi(a)]$ as
 divisors, or equivalently $[p^v](t)$ is a unit multiple of
 $\prod_a(t-x_a)$ in $D_v\psb{t}$.  It is also known that $x_a-x_b$ is
 invertible in $L$ whenever $a\neq b$ (because it is a unit multiple
 of $x_{a-b}$, which divides
 $\prod_{c\neq 0}x_c=\pm[p^v]'(0)=\pm p^v$).  Let the representations
 $V_i$ and the elements $c_{ik}\in C(E,G)$ be as in the proof of
 Proposition~\ref{prop-scheme}.  If $f\in\OmCh(G)$ and
 $f(V_i)=[a_1]+\ldots+[a_d]\in\N[\Tht(v)]$ then $\tCh(c_{ik})(f)$ is
 the $k$'th symmetric function in the variables
 $x_{a_1},\ldots,x_{a_d}$, and this characterises $\tCh$.

 We next show that $\tCh$ is surjective.  For any $u\in\OmCh(G)$ we
 define $\ep_u\:C(L,G)\xra{}L$ by $\ep_u(c)=\tCh(c)(u)$, and we put
 $I_u=\ker(\ep_u)\leq  C(L,G)$.  By the Chinese Remainder Theorem, it
 will suffice to show that $I_u+I_v=C(L,G)$ whenever $u\neq v$.  If
 $u\neq v$ we can choose $V\in R^+(G)$ such that
 $u(V)\neq v(V)\in\N[\Tht]$.  If $u(V)=\sum_am_a[a]$ and
 $v(V)=\sum_an_a[a]$ then we must have $m_b\neq n_b$ for some $b$, and
 without loss we may assume $m_b>n_b$.  Define
 \begin{align*}
  k_a &= \min(n_a,m_a) \\ 
  C   &= \sum_a k_a[a] \\
  A   &= \sum_a (n_a-k_a)[a] \\
  B   &= \sum_a (m_a-k_a)[a].
 \end{align*}
 We can write $A$ in the form $[a_1]+\ldots+[a_d]$ with $a_i\neq b$
 for all $i$.  We also write $f_A(t)=\prod_i(t-x_{a_i})$, so
 $f_A(x_b)$ is invertible in $L$.  On the other hand, the
 representation $V\in R_d^+(G)$ gives rise in a tautological way to a
 divisor $D_V\in\Div_d^+(\GG)(C(E,G))$ with equation
 $f_{D_V}(t)\in C(E,G)[t]$, say.  We have
 $\ep_uf_{D_V}(t)=f_{A+C}(t)=f_A(t)f_C(t)$ and
 $\ep_vf_{D_V}(t)=f_B(t)f_C(t)$.  The polynomial $f_C(t)$ is monic and
 thus is not a zero-divisor, so $f_A(t)=f_B(t)\pmod{I_u+I_v}$.  We
 evidently have $f_B(x_b)=0$ so $f_A(x_b)=0\pmod{I_u+I_v}$.  As
 $f_A(x_b)$ is invertible in $L$, we deduce that $I_u+I_v=1$ as
 required.

 Finally, we must show that the kernel of $\tCh$ is nilpotent.  This
 kernel is the intersection of the ideals $I_u$, so by well-known
 arguments it suffices to show that every prime ideal in $C(L,G)$
 contains $I_u$ for some $u$.  To see this, put $R=C(L,G)$ and let
 $\pri\leq R$ be a prime ideal.  If $D\in\Div_d^+(\GG)(C(E,G))$ is a
 divisor satisfying $\psi^{p^v}(D)=d[0]$, then Theorem~\ref{thm-Zdv}
 implies that $D\in\GG(v,d)^d/\Sg_d$ and thus that $D\leq d^2.\GG(v)$
 as divisors, or equivalently $f_D(t)$ divides $[p^v](t)^{d^2}$, which
 is a unit multiple in $D_v[t]$ of $\prod_{a\in\Tht(v)}(t-x_a)^{p^v}$.
 Now let $K$ be the field of fractions of $R/\pri$, and note that
 $x_a-x_b$ is invertible in $L$ and thus in $K$ when $a\neq b$.  As
 $K[t]$ is a unique factorisation domain, we see that
 $f_D(t)=\prod_a(t-x_a)^{m_a}$ in $K[t]$ for a unique system of
 integers $m_a$.  We define $w(D)=\sum_am_a[a]\in\Z[\Tht(v)]^+_d$.  In
 particular, if $V\in R_d^+(G)$ we can let $D_V$ be the tautologically
 associated divisor over $C(E,G)$ and put $u(V)=w(D_V)$.  One can
 check that this gives a homomorphism $u\:R^+(G)\xra{}\Z[\Tht]$ of
 $\Lm$-semirings, or in other words an element $u\in\OmCh(G)$.  From
 the construction it is automatic that $I_u\leq\pri$.
\end{proof}

\begin{example}
 If $G$ is Abelian, it is easy to see that
 $\OmCh(G)=\Hom(G^*,\Tht)\simeq\Hom(\Tht^*,G)=\Om(G)$ and that $\kp$
 is an isomorphism.
\end{example}
\begin{example}\label{eg-Sg-kp}
 Consider the symmetric group $G=\Sg_k$.  This acts in an obvious way
 on $\C^k$, and we call this representation $\pi$.  It is known 
 that $\pi$ generates $R(G)$ as a $\Lm$-ring.  Thus, an element
 $f\in\OmCh(G)$ is determined by the value $f(\pi)\in\Z[\Tht]_k^+$.

 As discussed in~\cite{sttu:rme}, the set $\Om(G)$ can be identified
 with the set of isomorphism classes of sets of order $k$ with an
 action of $\Tht^*$.  For any finite subgroup $A<\Tht$ we have a
 homomorphism $\Tht^*\xra{}A^*$ and thus an action of $\Tht^*$ on
 $A^*$.  Note that $0^*$ is just a single point with trivial action.
 We write $m.A^*$ for the disjoint union of $m$ copies of $A^*$.  If
 $T$ is a finite $\Tht^*$-set, then $T$ can be written in an
 essentially unique way in the form $\coprod_im_i.A_i^*$.

 If we write $[A]=\sum_{a\in A}[a]\in\Z[\Tht]$ then by working through
 the definitions we find that
 $\kp(\coprod_im_i.A_i^*)(\pi)=\sum_im_i[A_i]$, which effectively
 determines $\kp$.

 It is now easy to exhibit cases in which $\kp$ is not injective.  For
 example, suppose that $p=2$ and $n>1$ and $k=6$.  We can then find
 two distinct, nonzero elements $a,b\in\Tht(1)$ and put $c=a+b$.  Let
 $A$, $B$ and $C$ be the groups generated by $a$, $b$ and $c$
 respectively, and put $V=A+B=\{0,a,b,c\}$.  Then
 \[ \kp(V^*\amalg 2.0^*)(\pi)=
    \kp(A^*\amalg B^*\amalg C^*)(\pi)=  3[0] + [a] + [b] +[c],
 \]
 so $\kp$ is not injective.  In Section~\ref{sec-xspec} we will give
 examples where $\kp$ is not surjective.
\end{example}

\section{Calculating $\OmCh(G)$}\label{sec-OmChG}

In this section we define sets $\Om'(G)$ and $\Om''(G)$ which in some
cases may be easier to compute than $\Om(G)$ or $\OmCh(G)$, and we
define natural maps
\[ \Om(G)\;\arrow{e,A}\;
   \Om'(G)\;\arrow{e,V}\;
   \OmCh(G)\;\arrow{e,V}\; \Om''(G). 
\]

\begin{definition}
 Let $C$ be the set of conjugacy classes of elements of $p$-power
 order in $G$.  We let the multiplicative monoid $\Z$ act on $\Tht^*$
 in the obvious way, and on $C$ by $k.[g]=[g^k]$.  We say that two
 homomorphisms $u,v\:\Tht^*\xra{}G$ are pointwise conjugate if $u(a)$
 is conjugate to $v(a)$ for all $a\in\Tht^*$.  We recall the
 definitions of $\Om(G)$ and $\OmCh(G)$ and define new sets $\Om'(G)$
 and $\Om''(G)$ as follows:
 \begin{align*}
  \Om(G)   &= \Hom(\Tht^*,G)/\text{conjugacy} \\
  \Om'(G)  &= \Hom(\Tht^*,G)/\text{pointwise conjugacy} \\
  \Om''(G) &= \{\text{$\Z$-equivariant continuous maps }
                \Tht^* \xra{} C \} \\
  \OmCh(G) &= \{\text{ positive $\Lm$-ring homomorphisms }
                R(G) \xra{}\Z[\Tht] \}.
 \end{align*}
\end{definition}

\begin{proposition}
 There are natural maps as follows:
 \[ \Om(G)\;\arrow{e,A}\; 
    \Om'(G)\;\arrow{e,V}\;
    \OmCh(G)\;\arrow{e,V}\;\Om''(G).
 \]
\end{proposition}
\begin{proof}
 There are evident natural maps
 \[ \Om(G) \;\arrow{e,A}\; \Om'(G) \;\arrow{e,V}\; \Om''(G). \]
 We have also already constructed a map $\kp\:\Om(G)\xra{}\OmCh(G)$.
 If $u,v\:\Tht^*\xra{}G$ are pointwise-conjugate then the induced maps
 from class functions on $G$ to class functions on $\Tht^*$ are
 evidently the same, so the induced maps $R(G)\xra{}\Z[\Tht]$ are the
 same, so $\kp(u)=\kp(v)$.  This shows that $\kp$ factors through the
 projection $\Om(G)\xra{}\Om'(G)$.  

 We next define a map $\xi\:\OmCh(G)\xra{}\Om''(G)$.  Suppose that
 $u\in\OmCh(G)$ and $a\in\Tht^*=\Hom(\Tht,S^1)$.  Then $a$ extends in
 a natural way to give a $\C$-algebra map $\ha\:\C[\Tht]\xra{}\C$ and
 thus a ring map $(1_\C\ot u)\ha\:\C\ot R(G)\xra{}\C$.  Using the fact
 that $\C\ot R(G)$ is the set of $\C$-valued class functions on $G$,
 we see that $\Hom_{\C-\text{Alg}}(\C\ot R(G),\C)$ can be identified
 with the set of conjugacy classes in $G$.  Thus there exists $h\in G$
 (unique up to conjugation) such that $(1\ot u)(\ha(V))=\chi_V(h)$ for
 all $V\in R(G)$.  We can choose $m$ so that $u(V)\in\Z[\Tht(m)]$ for
 all $V$, and then we have
 $\chi_V(h^{p^m})=\chi_{\psi^{p^m}V}(h)=\chi_{\dim(V)}(h)=\dim(V)$ for
 all $V$, so $h^{p^m}=1$.  This means that the conjugacy class $[h]$
 lies in $C$, so we can define $\xi(u)(a)=[h]$.  We leave it to the
 reader to check that this gives a map $\xi\:\OmCh(G)\xra{}\Om''(G)$
 as claimed.  The maps $\ha\:\Z[\Tht]\xra{}\C$ (as $a$ runs over
 $\Tht^*$) are jointly injective, an it follows that $\xi$ is
 injective.  One can also check that the composite
 $\Om'(G)\xra{}\OmCh(G)\xra{}\Om''(G)$ is just the obvious inclusion,
 which implies that the map $\Om'(G)\xra{}\OmCh(G)$ is injective.
\end{proof}

\section{Special divisors}

In this section we study ``special'' divisors, which are related to
the special unitary group in the same way that arbitrary divisors are
related to the full unitary group.
\begin{definition}\label{defn-SDiv}
 A divisor $D\in\Div_d^+(\GG)$ is \emph{special} if $\lm^dD=[0]$.  We
 write $\SDiv_d^+(\GG)$ for the scheme of special divisors.
\end{definition}

\begin{proposition}\label{prop-SDiv}
 We have $\O_{\SDiv_d^+(\GG)}=\OX\psb{c_2,\ldots,c_d}$.  In the
 topological situation this can be identified with $E^0BSU(d)$.  
\end{proposition}
\begin{proof}
 Put $A=\O_{\GG^d}=\OX\psb{x_1,\ldots,x_d}$ and
 $A'=\O_{\Div_d^+(\GG)}=A^{\Sg_d}=\OX\psb{c_1,\ldots,c_d}$.  Here
 $c_i$ is the $i$'th elementary symmetric function, and in particular
 $c_1=\sum_ix_i$.  Put $c'_1=\sum^F_ix_i\in A'$.  If we regard $\lm^d$
 as a map $\Div_d^+(\GG)\xra{}\Div_1^+(\GG)=\GG$ then
 $c'_1=x\circ\lm^d$, so we see that $\O_{\HH_3}=A/c'_1$ and
 $\O_{\SDiv_d^+(\GG)}=A'/c'_1$.  Next, observe that the inclusion
 $A'\xra{}A$ induces an inclusion
 $A'/(c_1^2,c_2,\ldots,c_d)\xra{}A/(x_1,\ldots,x_d)^2$.  We have
 $c'_1=c_1\pmod{(x_1,\ldots,x_d)^2}$ so
 $c'_1=c_1\pmod{c_1^2,c_2,\ldots,c_d}$.  It follows easily that
 $A'=\OX\psb{c'_1,c_2,\ldots,c_d}$ and thus that
 $A'/c'_1=\OX\psb{c_2,\ldots,c_d}$.
\end{proof}

We next put $\HH_d=\ker(\GG^d\xra{+}\GG)$.  If we let
$q\:\GG^d\xra{}\Div_d^+(\GG)$ be the usual projection (which is finite
and faithfully flat, with degree $d!$) then
$\HH_d=q^{-1}\SDiv_d^+(\GG)$.  It follows that the map
$q\:\HH_d\xra{}\SDiv_d^+(\GG)$ is also finite and faithfully flat,
with the same degree.  It clearly factors through
$\HH_d/\Sg_d:=\spf(\O_{\HH_d}^{\Sg_d})$, and one would like the
induced map $q\:\HH_d/\Sg_d\xra{}\SDiv_d^+(\GG)$ to be an isomorphism.
However, quotient constructions in algebraic geometry are never as
simple as one would like, and we do not know whether this is true in
general; certainly it becomes false if we remove our assumption that
$\GG$ has finite height.  For example, consider the case where $\GG$
is the additive group over $\F_2$ and $d=2$; then $2a=0$ for all
$a\in\GG$ so $\Sg_2$ acts trivially on $\HH_2=\{(a,-a)\st a\in\GG\}$
so the map $q\:\HH_d/\Sg_d\xra{}\SDiv_d^+(\GG)$ has degree two.
However, we do have the following partial result.
\begin{proposition}\label{prop-HdSgd}
 If $d$ is invertible in $\OX$ then $\SDiv_d^+(\GG)=\HH_d/\Sg_d$.
\end{proposition}
\begin{proof}
 As $d$ is invertible in $\OX$, multiplication by $d$ is an
 automorphism of $\GG$.  Define maps
 $\GG\tm\HH_d\xra{f}\GG^d\xra{g}\GG$ by
 $f(a,b_1,\ldots,b_d)=(a+b_1,\ldots,a+b_d)$ and
 $g(b_1,\ldots,b_n)=\sum_ib_i/d$, and then define
 $h\:\GG^d\xra{}\GG\tm\HH_d$ by
 $h(\ub)=(g(\ub),b_1-g(\ub),\ldots,b_d-g(\ub))$.  Clearly $h$ is
 inverse to $f$, so $f$ is an isomorphism, giving an isomorphism
 $\O_{\GG^d}=\O_G\hot\O_{\HH_d}=\O_{\HH_d}\psb{x}$ of rings.  If we
 let $\Sg_d$ act trivially on $\GG$ then everything is equivariant, so
 we have $\O_{\GG^d}^{\Sg_d}=\O_{\HH_d}^{\Sg_d}\psb{x}$, so
 $\Div_d^+(\GG)=\GG^d/\Sg_d=\GG\tm(\HH_d/\Sg_d)$.
\end{proof}

\section{The group $\Sg_4$}\label{sec-Sigma-four}

We now consider the case where $G=\Sg_4$ and $E$ is the $2$-periodic
Morava $E$-theory spectrum of height $2$ at the prime $2$.  We shall
show that the map $C(E,G)\xra{}E^0BG$ is an isomorphism.  To be more
explicit, we need to name some representations.  Note that $\Sg_4$
acts on $\C^4$ with a one-dimensional fixed subspace; we let $\rho$ be
the representation of $G$ on the quotient space.  We also write $\ep$
for the sign representation.  We let $K=E/I_2$ denote the $2$-periodic
Morava $K$-theory spectrum.  
\begin{theorem}\label{thm-Sg-four}
 Let $c_2,c_3\in E^0B\Sg_4$ be the Chern classes of the representation
 $\ep\rho$, and let $w$ be the Euler class of $\ep$.  Then
 $C(E,\Sg_4)=E^0B\Sg_4$, and this is a free module of rank $17$ over
 $E^0$, with the following monomials as a basis:
 \[ \setlength{\arraycolsep}{2em}
 \begin{array}{ccccc}
  1   &     c_2 &     c_2^2 & c_2^3 &     c_3 \\
  w   & w   c_2 & w   c_2^2 &       & w   c_3 \\
  w^2 & w^2 c_2 & w^2 c_2^2 &       & w^2 c_3 \\
  w^3 & w^3 c_2 & w^3 c_2^2 &       & w^3 c_3.
 \end{array} \]
 Moreover, we have 
 \[ C(K,\Sg_4)=K^0B\Sg_4=C(E,\Sg_4)/I_2=K^0[w,c_2,c_3]/J, \]
 where $J$ is generated by the following elements:
 \begin{align*}
  & w^4 \;,\; c_3^2 \;,\; c_2c_3 \;,            \\
  & c_2^4 + w^2 c_2^3 + w c_2^2 + w^2 c_3,      \\
  & w c_2^3 + w^2 c_2 + w c_3.  
 \end{align*}
\end{theorem}
We will prove this in a number of stages.  In
Section~\ref{subsec-R-Sg-four} we assemble the facts that we need
about the representations of $\Sg_4$, and in
Section~\ref{subsec-HKR-Sg-four} we deduce that the map
$\kp\:\Om(\Sg_4)\xra{}\OmCh(\Sg_4)$ is a bijection.  We then recall
some formulae for the relevant formal group law, and in
Section~\ref{subsec-SDiv-three} we use them to analyse the structure
of an auxiliary scheme denoted $\SDiv_3^+(\GG_0)^C$.  This allows us
to complete our determination of $C(K,\Sg_4)$ in
Section~\ref{subsec-C-K-Sg-four}, with the help of some theory of
Gr\"obner bases.  We find in particular that $C(K,\Sg_4)$ is a
Gorenstein ring, which enables us to use the inner products defined
in~\cite{st:kld} to show that the map
$\tht\:C(K,\Sg_4)\xra{}K^0B\Sg_4$ is injective; this is explained in
Section~\ref{subsec-transfer}.  We know from~\cite{hokura:ggc} that
$K(n)^*B\Sg_4$ is concentrated in even degrees, and it follows that
$E^0B\Sg_4$ is a free module over $E^0$ of rank $|\Om(\Sg_4)|=17$;
see~\cite{st:msg} for more details.  In Section~\ref{subsec-proof} we
combine these various facts to prove the theorem.

\subsection{Representation theory}\label{subsec-R-Sg-four}

Our first task is to understand the structure of $R(\Sg_4)$.  We have
already defined the characters $\ep$ and $\rho$.  It is a standard
calculation that there is another irreducible character $\sg$ of
dimension $2$ such that the character table is as follows:
\[ \setlength{\arraycolsep}{2em}
   \renewcommand{\arraystretch}{1.5}
   \begin{array}{|c|c|ccccc|}
   \hline
   \text{class}&\text{size}& 1 & \ep & \sg  & \rho& \ep\rho \\ \hline
   1^4         & 1         & 1 & 1   & 2    & 3   & 3      \\ 
   1^22        & 6         & 1 &-1   & 0    & 1   &-1      \\ 
   2^2         & 3         & 1 & 1   & 2    &-1   &-1      \\ 
   13          & 8         & 1 & 1   &-1    & 0   & 0      \\ 
   4           & 6         & 1 &-1   & 0    &-1   & 1      \\ \hline
   \end{array}
\]

The ring structure, Adams operations and $\Lm$-operations are
described in the following table.
\[ \setlength{\arraycolsep}{2em}
   \renewcommand{\arraystretch}{1.5}
   \begin{array}{lll}
    \ep^2  =1                         &
     \psi^k(\ep) =\ep^k                &
      \lm^2(\sg) = \ep                  \\
    \ep\sg =\sg                       &
     \psi^2(\sg) =1-\ep+\sg            &   
      \lm^2(\rho) = \ep\rho             \\
    \sg^2  =1+\ep+\sg                 &
     \psi^3(\sg) =1+\ep                &   
      \lm^3(\rho) = \ep                 \\
    \sg\rho=\rho+\ep\rho              &
     \psi^2(\rho)=1+\sg+\rho-\ep\rho   & \\
    \rho^2 =1+\sg+\rho+\ep\rho        &
     \psi^3(\rho)=1+\ep-\sg+\rho.      &
   \end{array}
\]
(The first two columns are easily checked by looking at the
characters, and the last column follows using the standard formulae
relating Adams operations to $\Lm$-operations.)

Let $P$ be a Sylow $2$-subgroup (a dihedral group of order $8$) and
$I$ be the kernel of the restriction map $R(\Sg_4)\xra{}R(P)$; one
checks that $I=(\sg-1-\ep)$.  Put $\tau=\ep\rho\in R_3^+(\Sg_4)$; one
checks that $\lm^k(\tau)=\ep^k\lm^k(\rho)$ and so $\lm^2(\tau)=\tau$
and $\lm^3\tau=1$.  We have
\[ R(\Sg_4)/I = \Z\{1,\ep,\tau,\ep\tau\} =
   \Z[\ep,\tau]/(\ep^2-1,\tau^2-1-(1+\ep)(1+\tau)).
\]
The operation $\psi^k$ acts as the identity on this ring when $k$ is
odd, and we have
\begin{align*}
 \psi^2(\ep)  &= 1 \\
 \psi^2(\tau) &= 2+\ep+\ep\tau-\tau.
\end{align*}

\begin{proposition}\label{prop-OmCh-Sg-four}
 The set $\OmCh(\Sg_4)$ can be identified with the set of pairs
 $(d,u)\in\Tht(1)\tm\Z[\Tht]_3^+$ such that 
 \begin{align*}
  2d           &= 0                 \\
  \lm^3(u)     &= [0]               \\
  \psi^{-1}(u) &= u                 \\
  \psi^2(u)+u  &= 2[0] + [d] + [d]u .
 \end{align*}
 Similarly, $\XCh(\Sg_4)$ can be identified with the scheme of pairs
 $(d,D)\in\GG(1)\tm\Div_3^+(\GG)$  such that
 \begin{align*}
  2d           &= 0                 \\
  \lm^3(D)     &= [0]               \\
  \psi^{-1}(D) &= D                 \\
  \psi^2(D)+D  &= 2[0] + [d] + [d]D.
 \end{align*}
\end{proposition}
\begin{proof}
 Given a positive homomorphism $f\:R(\Sg_4)\xra{}\Z[\Tht]$, let
 $d\in\Tht$ be the element such that $f(\ep)=[d]$ and put $u=f(\tau)$.
 We know from Remark~\ref{rem-OmCh-i} that $f(I)=0$ and it follows
 easily from our description of $R(\Sg_4)/I$ that $d$ and $u$ have the
 properties listed.  Conversely, given $d$ and $u$ as described, we
 can define a homomorphism $f\:R(\Sg_4)/I\xra{}\Z[\Tht]$ of additive
 groups by
 \begin{align*}
  f(1)       &= [0]     \\
  f(\ep)     &= [d]     \\
  f(\tau)    &= D       \\
  f(\ep\tau) &= [d]D.
 \end{align*}
 It is straightforward to check that this gives a homomorphism of
 $\Lm$-rings, and that these constructions give the required
 bijection.  The argument for $\XCh(\Sg_4)$ is essentially the same.  
\end{proof}

\subsection{Generalised character theory}\label{subsec-HKR-Sg-four}

We next work out the generalised character theory (as recalled in
Section~\ref{sec-gen-char}) of $\Sg_4$.  The set $\Om(\Sg_4)$ can be
described in terms of $\Tht^*$-sets as in Example~\ref{eg-Sg-kp}.  We
can thus write $\Om(\Sg_4)$ as the disjoint union
$\Phi_0\amalg\ldots\amalg\Phi_4$, where
\begin{itemize}
 \item $\Phi_0$ consists of the set
  $4.0^*:=0^*\amalg 0^*\amalg 0^*\amalg 0^*$.
 \item $\Phi_1$ consists of the sets $2.0^*\amalg A^*$, where
  $A\simeq \Z/2$ (so $|\Phi_1|=2^n-1$).
 \item $\Phi_2$ consists of the sets $A^*\amalg B^*$, where
  $A\simeq B\simeq\Z/2$, and $A$ may be equal to $B$.  We have
  $|\Phi_2|=\half|\Phi_1|(|\Phi_1|+1)=2^{n-1}(2^n-1)$.
 \item $\Phi_3$ consists of the sets $A^*$ where $A\simeq(\Z/2)^2$.
  We have $|\Phi_3|=(2^n-1)(2^{n-1}-1)/3$ (by counting the number of
  linearly independent pairs in $(\Z/2)^n$ and dividing by
  $|GL_2(\Z/2)|=6$).
 \item $\Phi_4$ consists of the sets $A^*$ where $A\simeq\Z/4$.  There
  are $2^{2n}-2^n$ points in $\Tht$ of order exactly $4$, and each
  subgroup in $\Phi_4$ contains precisely two of these, so 
  $|\Phi_4|=(2^{2n}-2^n)/2=2^{n-1}(2^n-1)$.
\end{itemize}

\begin{proposition}\label{prop-HKR-Sg-four}
 The map $\kp\:\Om(\Sg_4)\xra{}\OmCh(\Sg_4)$ is a bijection.
\end{proposition}
\begin{proof}
 Define 
 \[ \Psi_i=\kp(\Phi_i)\sse\OmCh(\Sg_4). \]
 Recall from Example~\ref{eg-Sg-kp} that
 $\kp(\coprod_im_i.A_i^*)(\pi)=\sum_im_i[A_i]$.  An easy case-by-case
 check shows that the sets $\Psi_i$ are disjoint and that the maps
 $\kp\:\Phi_i\xra{}\Psi_i$ are bijections.  It will thus be enough to
 show that the union of the sets $\Psi_i$ is the whole of
 $\OmCh(\Sg_4)$.

 Suppose we have an element $f\in\OmCh(\Sg_4)$, with $f(\ep)=[d]$ and
 $f(\tau)=u=[a]+[b]+[c]$ say.  Let $A$, $B$ and $C$ be the cyclic
 subgroups generated by $a$, $b$ and $c$ respectively.  Put
 $v=f(\pi)=[0]+[d]u=[0]+[a+d]+[b+d]+[c+d]$, and recall that this
 determines $f$, because $\pi$ generates $R(\Sg_4)$ as a $\Lm$-ring.  As
 $\psi^4\tau=3[0]$ we have $4a=4b=4c=0$.  By
 Proposition~\ref{prop-OmCh-Sg-four}, we have 
 \begin{align*}
  2d          &= 0 \\
  a+b+c       &= 0 \\
  [a]+[b]+[c] &= [-a]+[-b]+[-c] \\
  [2a]+[2b]+[2c]+[a]+[b]+[c] &= 2[0]+[d]+[a+d]+[b+d]+[c+d].
 \end{align*}
 Suppose that $\psi^2(u)\neq 3[0]$.  Without loss of generality we may
 assume that $2a\neq 0$ so $-a\neq a$.  The third equation implies
 that $-a\in\{b,c\}$, so we may assume that $-a=b$.  As $a+b+c=0$ we
 must have $c=0$.  Recall also that $4a=0$ so $2a=-2a$.  Putting all
 this in the last equation and cancelling $2[0]$ gives
 \[ 2[2a]+[a]+[-a] = 2[d] + [d+a] + [d-a]. \]
 Note that $2a$ and $d$ have order $2$, but $a$, $-a$, $d+a$ and $d-a$
 do not.  It follows that we must have $2a=d$ and thus
 $v=[0]+[a]+[2a]+[3a]$.  We conclude that $f=\kp(A^*)\in\Psi_4$.

 We may thus assume that $\psi^2(u)=3[0]$, so $2a=2b=2c=0$.  Suppose
 that $d=0$.  As $a+b+c=0$ we see that $D:=\{0,a,b,c\}$ is a subgroup
 of $\Tht$, of order $2^e$ say (so $e\in\{0,1,2\}$).  This implies
 that $v=2^{2-e}[A]=\kp(2^{2-e}.D^*)(\pi)$, so
 $f=\kp(2^{2-e}.D^*)\in\Psi_0\amalg\Psi_2\amalg\Psi_3$.
 
 We may thus assume that $2a=2b=2c=2d=0$ and $d\neq 0$.  The equation
 $\psi^2(u)+u=2[0]+[d]+[d]u$ then reduces to
 \[ [0] + [a] + [b] + [c] = [d] + [a+d] + [b+d] + [c+d]. \]
 It follows that $d\in\{a,b,c\}$ and without loss we may assume that
 $d=a$.  Note that $c=a+b=d+b$ (because $a+b+c=0$).  If $b=0$ this
 gives $c=a=d$ so $v=3[0]+[a]$, so $f=\kp(2.0^*\amalg A^*)$.  The same
 argument works if $c=0$, so we reduce to the case where $a=d\neq 0$
 and $b$ and $c$ are also nonzero.  We then have
 \[ v=[0]+[a+d]+[b+d]+[c+d]= 2[0]+[c]+[b]=[B]+[C], \]
 so $f=\kp(B^*\amalg C^*)\in\Phi_2$.
\end{proof}

\subsection{The formal group law}\label{subsec-FGL}

Let $\GG$ be the formal group associated to $E$, and let $\GG_0$ be
its restriction to the special fibre $X_0\subset X$, or equivalently
the formal group associated to $K$.  This has a standard coordinate
giving rise to a formal group law $F$ over $\O_{X_0}=K^0=\F_4$, which
is in fact defined over $\F_2$.  We will need the following formulae:
\begin{align*}
 [2](x)  &= x^4                                                 \\
 [-1](x) &= x + x^4 + x^{10} + x^{16} + x^{22} \pmod{x^{32}}    \\
 x +_F y &= x + y + x^2 y^2 \pmod{x^4 y^4}.
\end{align*}
The first of these is well-known and the second can be proved by
straightforward computation; for the third, one can adapt the method
of~\cite[Section 15]{st:fsf} to the case $p=2$.

\subsection{The scheme $\SDiv_3^+(\GG_0)^C$}\label{subsec-SDiv-three}

Let $C$ be the group (of order $2$) generated by $\psi^{-1}$, so 
\[ \SDiv_3^+(\GG_0)^C =
   \{D \in \Div_3^+(\GG_0)\st
     \lm^3D=[0] \text{ and } \psi^{-1}D=D \}.
\]

We have seen that $\O_{\SDiv_3^+(\GG_0)}=\F_4\psb{c_2,c_3}$; our next
task is to determine the quotient ring $\O_{\SDiv_3^+(\GG_0)^C}$.
\begin{proposition}\label{prop-C-invariant}
 We have 
 \[ \O_{\SDiv_3^+(\GG_0)^C}=\F_4\psb{c_2,c_3}/(c_2c_3,c_3^2)=
    \F_4\psb{c_2}\oplus\F_4.c_3.
 \]
\end{proposition}
\begin{proof}
 Put $A=\F_4\psb{x,y,z}$ and
 \begin{align*}
  d   &= x+_Fy+_Fz      \\
  c_1 &= x+y+z          \\
  c_2 &= xy + yz + zx   \\
  c_3 &= xyz.
 \end{align*}
 Put $A'=A^{\Sg_3}=\F_4\psb{d,c_2,c_3}$, so $A$ is free of rank $6$
 over $A'$.  Put $B'=A'/dA'$ and
 \[ B=A/dA=B'\ot_{A'}A=\F_4\psb{c_2,c_3}=\O_{\SDiv_3^+(\GG_0)}. \]
 For any element $u\in A$ we write $\ov{u}=(\psi^{-1})^*(u)$, so
 $u\mapsto\ov{u}$ is a ring map and $\ov{u}=[-1](u)$ for
 $u\in\{x,y,z,d\}$.  Put $C=B/(\cb_2-c_2,\cb_3-c_3)B$ and
 $C'=B'/(\cb_2-c_2,\cb_3-c_3)B'=\O_{\SDiv_3^+(\GG_0)^C}$.  The claim
 is that the ideal in $C'$ generated by $c_3$ is free of rank one over
 $\F_4$, and that $C'/c_3C'=\F_4\psb{c_2}$, so that
 $C'=\F_4\psb{c_2}\oplus\F_4.c_3$.

 We will think of $\Div_2^+(\GG_0)$ as being embedded in
 $\Div_3^+(\GG_0)$ by the map $D\mapsto D+[0]$, so
 \[ \O_{\Div_2^+(\GG_0)}=\O_{\Div_3^+(\GG_0)}/c_3=\F_4\psb{d,c_2}. \]
 There is a faithfully flat map $\GG_0\xra{}\SDiv_2^+(\GG_0)$ sending
 $a$ to $[a]+[-a]$, and clearly $\psi^{-1}([a]+[-a])=[a]+[-a]$ so
 $\SDiv_2^+(\GG_0)\sse\SDiv_3^+(\GG_0)^C$.  It follows that
 $\Div_2(\GG_0)\cap\SDiv_3^+(\GG_0)^C=\SDiv_2^+(\GG_0)$, and thus that
 $C'/c_3C'=\F_4\psb{c_2}$ as claimed.

 This implies that we must have $\cb_2-c_2=c_3r_2$ and
 $\cb_3-c_3=c_3r_3$ for some $r_2,r_3\in B'$.

 Now work in $B/(x,y,z)^7$.  We have
 \begin{align*}
  z   &= \xb +_F \yb = x + y + x^4 + x^2 y^2 + y^4      \\
  \xb &= x + x^4                                        \\
  \yb &= y + y^4                                        \\
  \zb &= x + y + x^2 y^2                                \\
  c_2 &= x^2 + x y + y^2 +
         x^5 + x^4 y + x^3 y^2 + x^2 y^3 + x y^4 + y^5  \\
  c_3 &= x^2 y + x y^2 + x^5 y + x^3 y^3 + x y^5        \\
  \cb_2-c_2 &= c_2c_3 = x^4 y + x y^4                   \\
  \cb_3-c_3 &= c_3^2 = x^4 y^2 + x^2 y^4.
 \end{align*}
 We also find that the ideal $c_3.(c_2,c_3)^2$ maps to zero in this
 ring.  Using this, we find that $r_2=c_2\pmod{(c_2,c_3)^2}$ and
 $r_3=c_3\pmod{(c_2,c_3)^2}$, so $B'=\F_4\psb{r_2,r_3}$ and
 $B'/(r_2,r_3)=\F_4$.  It follows that
 $\O_{\SDiv_3^+(\GG_0)^C}=B'/(r_2c_3,r_3c_3)=B'/(c_2c_3,c_3^2)3$ as
 claimed.
\end{proof}

Now put $Y=\{D\in\SDiv_3^+(\GG_0)^C\st\psi^4D=3[0]\}$, and let
$U\subset\GG_0$ be the divisor $32[0]$.  We know that
$c_k(\psi^4D)=c_k(D)^{16}$ so
\[ \OY=\O_{\SDiv_3^+(\GG_0)^C}/(c_1^{16},c_2^{16},c_3^{16})=
      \F_4[c_2,c_3]/(c_2^{16},c_3^2,c_2c_3). 
\]
We can also study $Y$ using the maps
\[ \al\:U\xra{}Y \hspace{4em} \al(a)=[a]+[-a]+[0] \]
\[ \bt\:\GG_0(1)^2\xra{}Y \hspace{4em} \bt(a,b)=[a]+[b]+[a+b]. \]

The map $\al$ gives a ring map $\al^*\:\OY\xra{}\F[x]/x^{32}$, with
\begin{align*}
 \al^*(c_1) &= x+\xb = x^4 + x^{10} + x^{16} + x^{22}           \\
 \al^*(c_2) &= x \xb = x^2 + x^5 + x^{11} + x^{17} + x^{23}     \\
 \al^*(c_3) &= 0. 
\end{align*}
The map $\bt$ gives a ring map $\bt^*\:\OY\xra{}\F[x,y]/(x^4,y^4)$.
If we put $z=x+_Fy=x+y+x^2y^2$ then 
\begin{align*}
 \bt^*(c_1) &= x + y + z = x^2y^2                               \\
 \bt^*(c_2) &= xy + yz + zx =  x^2 + xy + y^2 + x^2y^2(x+y)     \\
 \bt^*(c_3) &= xyz = xy(x+y) + x^3y^3.
\end{align*}

\begin{proposition}\label{prop-al-bt}
 The maps $\al^*$ and $\bt^*$ are jointly injective (in other words,
 $\ker(\al^*)\cap\ker(\bt^*)=0$).  Moreover, we have $c_1=c_2^2+c_2^8$
 in $\OY$.
\end{proposition}
\begin{proof}
 Recall that $\OY=\F_4[c_2,c_3]/(c_2^{16},c_3^2,c_2c_3)$, so
 $\{c_2^i\st 0\leq i<16\}\amalg\{c_3\}$ is a basis for $\OY$ over
 $\F_4$.  As $\al^*(c_2)=x^2\pmod{x^3}$, it is easy to see that
 $\al^*(c_2^i)=x^{2i}\pmod{x^{2i+1}}$ and that these elements are
 linearly independent in $\O_U=\F_4[x]/x^{32}$.  Moreover, we have
 \begin{align*}
  \bt^*(1)      &= 1                                    \\
  \bt^*(c_2)    &= x^2 + xy + y^2 + x^2y^2(x+y)         \\
  \bt^*(c_2^2)  &= x^2y^2                               \\
  \bt^*(c_2^3)  &= x^3y^3                               \\
  \bt^*(c_2^i)  &= 0 \hspace{5em} \text{ for } i>3.
 \end{align*}
 It is easy to check that $\bt^*(c_3)$ does not lie in the span of
 these elements, and to deduce that $\al^*$ and $\bt^*$ are jointly
 injective as claimed.  Thus, to show that $c_1=c_2^2+c_2^8$ we need
 only check that $\al^*(c_1)=\al^*(c_2^2+c_2^8)$ and
 $\bt^*(c_1)=\bt^*(c_2^2+c_2^8)$, which is a straightforward
 computation. 
\end{proof}

\subsection{The ring $C(K,\Sg_4)$}\label{subsec-C-K-Sg-four}

Consider a pair $(d,D)\in\GG_0(1)\tm Y$.  This gives us a divisor
$[d]D\in\Div_3^+(\GG_0)$ defined over the ring
\[ \O_{\GG_0(1)}\ot\OY=\F_4[w,c_2,c_3]/(w^4,c_2^{16},c_3^2,c_2c_3). \]
Here of course $c_2$ and $c_3$ are the usual invariants of the divisor
$D$, but the divisor $[d]D$ also has invariants $c_k([d]D)$ lying in
the above ring.  In order to apply the description of $\XCh(\Sg_4)$ in
Proposition~\ref{prop-OmCh-Sg-four}, we will need to understand these
invariants.  

\begin{proposition}\label{prop-cdD}
 We have
 \begin{align*}
  c_1([d]D) &= c_2^2 + c_2^8 + w + c_2^4 w^2                    \\
  c_2([d]D) &= c_2 + (1 + c_2^3 + c_2^9 + c_3) w^2              \\
  c_3([d]D) &= c_3 + c_2 w + (c_2^2 + c_2^8) w^2 + 
               (1 + c_3 + c_2^3 + c_2^9) w^3.
 \end{align*}
\end{proposition}
\begin{proof}
 First recall that $u+_Fv=u+v+u^2v^2\pmod{u^4v^4}$ and $w^4=0$ so
 $w+_Fv=w+v+w^2v^2$ for any $v$.

 Next note that $[d]\al(a)=[d]([a]+[-a]+[0])=[d+a]+[d-a]+[d]$, so
 $\al^*(c_k([d]D))=c_k([d+a]+[d-a]+[d])$ is the $k$'th elementary
 symmetric function of $\{w+_Fx,w+_F\xb,w\}$, for example
 \begin{align*}
  \al^*c_1([d]D) &= w+(w+x+w^2x^2)+(w+\xb+w^2\xb^2)             \\
   &= x^4 + x^{10} + x^{16} + x^{22} + w + x^8 w^2 + x^{20} w^2 \\
   &= \al^*(c_2^2 + c_2^8) + w + \al^*(c_2^4) w^2. 
 \end{align*}
 By similar computations, our other two equations also become true
 when we apply $\al^*$.

 In the same way, we have $[d]\bt(a,b)=[d+a]+[d+b]+[d+a+b]$, so
 $\bt^*c_k([d]D)$ is the $k$'th elementary symmetric function of the
 list $\{w +_F x, w +_F y, w +_F x +_F y\}$, or equivalently the list 
 \[ \{w + x + w^2x^2,
      w + y + w^2y^2,
      w + x + y + w^2x^2 + w^2y^2 + x^2y^2\}.
 \]
 We thus have
 \begin{align*}
  \bt^*c_3([d]D) &= 
   (w+x+w^2x^2)(w+y+w^2y^2)(w+x+y+w^2x^2+w^2y^2+x^2y^2) \\
  &= (x^2 y + x y^2 + x^3 y^3) +
     (x^2 + x y + y^2 + x^3 y^2 + x^2 y^3) w +          \\
  &\qquad x^2 y^2 w^2 +  (1+ x^2 y + x y^2) w^3         \\
  &= \bt^*(c_3) + \bt^*(c_2) w + \bt^*(c_2^2 + c_2^8) w^2 + 
     \bt^*(1 + c_3 + c_2^3 + c_2^9) w^3.
 \end{align*}
 By similar computations, our other two equations also become true
 when we apply $\bt^*$.  As $\al^*$ and $\bt^*$ are jointly injective,
 it follows that our equations hold in $\O_{\GG_0(1)\tm Y}$ as claimed.
\end{proof}

\begin{proposition}\label{prop-K-Sg-four}
 Let $J$ be the ideal in $\F_4[w,c_2,c_3]$ generated by the elements 
 \begin{align*}
  & w^4 \;,\; c_3^2 \;,\; c_2c_3 \;,            \\
  & c_2^4 + w^2 c_2^3 + w c_2^2 + w^2 c_3,      \\
  & w c_2^3 + w^2 c_2 + w c_3.  
 \end{align*}
 Then $C(K,\Sg_4)=\F_4[w,c_2,c_3]/J$.  Moreover, the following
 monomials form a basis for this ring over $\F_4$, so it has dimension
 $17$.
 \[ \setlength{\arraycolsep}{2em}
 \begin{array}{ccccc}
  1   &     c_2 &     c_2^2 & c_2^3 &     c_3 \\
  w   & w   c_2 & w   c_2^2 &       & w   c_3 \\
  w^2 & w^2 c_2 & w^2 c_2^2 &       & w^2 c_3 \\
  w^3 & w^3 c_2 & w^3 c_2^2 &       & w^3 c_3
 \end{array} \]
\end{proposition}
\begin{proof}
 Proposition~\ref{prop-OmCh-Sg-four} is equivalent to the statement
 that 
 \[ X(\Sg_4) = \{(d,D)\in\GG_0(1)\tm Y\st D+\psi^2(D)=2[0]+[d]+[d]D\}. \]
 This means that $C(K,\Sg_4)$ is the largest quotient of
 $\O_{\GG_0(1)\tm Y}$ over which we have $g(t)=0$, where
 \[ g(t)=f_D(t)f_{\psi^2D}(t)-t^2(t+w)f_{[d]D}(t). \]
 Here we write $f_D(t)=t^3+c_1(D)t^2+c_2(D)t+c_3(D)$ and similarly for
 our other divisors.  As usual we write $c_k$ for $c_k(D)$, and we
 recall from Proposition~\ref{prop-al-bt} that $c_1=c_2^2+c_2^8$.  We
 also recall that $c_k(\psi^2D)=c_k(D)^4$, so that 
 \[ f_{\psi^2D}(t) = t^3 + c_1^4 t^2 + c_2^4 t + c_3^4 =
                     t^3 + c_2^8 t^2 + c_2^4 t.
 \]
 The polynomial $f_{[d]D}(t)=\sum_{k=0}^3c_k([d]D)t^{3-k}$ can be read
 off from Proposition~\ref{prop-cdD}.  Putting all this together and
 expanding it out, we find that $g(t)=\sum_{k=1}^4r_kt^{6-k}$, where
 \begin{align*}
  r_1 &= c_2^8 + c_2^4 w^2 \\
  r_2 &= c_2^4 w^3 + (c_2^9 + c_2^3 + c_3) w^2 + 
         (c_2^8 + c_2^2) w + (c_2^{10} + c_2^4) \\
  r_3 &= (c_2^8 + c_2^2) w^2 + (c_2^{12} + c_2^9 + c_2^6) \\
  r_4 &= (c_2^8 + c_2^2) w^3 + c_2 w^2 + c_3 w + c_2^5.
 \end{align*}
 We thus have
 \[ C(K,\Sg_4)=
    \F_4[w,c_2,c_3]/(w^4,c_2^{16},c_3^2,c_2c_3,r_1,r_2,r_3,r_4).
 \]
 As $1+c_2^6+c_2^{12}+w^3$ is invertible, we can replace $r_2$ by
 \[ r'_2:=(1+c_2^6+c_2^{12}+w^3)r_2=
    c_2^4 + w^2 c_2^3 + w c_2^2 + w^2 c_3,
 \]
 which is one of the relations in the statement of the theorem.  As
 $w^4=0$ we have $(r'_2)^2=r_1$ and $(r'_2)^4=c_2^{16}$, so the
 relations $r_1$ and $c_2^{16}$ are redundant.  Similarly, we can
 replace $r_4$ by the relation 
 \[ r'_4 := r_4 + (c_2 + w^2 + c2^4 w^3) r'_2 =
     w c_2^3 + w^2 c_2 + w c_3,
 \]
 which is another of the relations in the statement of the theorem.
 One can check that 
 \[ r_3 = 
    (1 + c_2^3 + c_2^6) 
     (c_2^2 ( 1 + (1 + c_2^6)(w^3 + c_2^2 w^2)) r'_2 + c_2 r'_4),
 \]
 so $r_3$ is redundant.  We deduce that 
 \[ C(K,\Sg_4) = \F_4[w,c_2,c_3]/(w^4,c_3^2,c_2c_3,r'_4,r'_2) \]
 as claimed.

 We next show that the $17$ monomials listed form a basis for this
 quotient ring.  We order the set of monomials in $w$, $c_2$ and $c_3$
 by saying that $c_2^ic_3^jw^k<c_2^{i'}c_3^{j'}w^{k'}$ iff $i<i'$ or
 ($i=i'$ and $j<j'$) or ($i=i'$ and $j=j'$ and $k<k'$).  We claim that
 our relations form a Gr\"obner basis for $J$ with respect to this
 ordering.  We first recall briefly what this means.  The list of
 leading terms of our relations is $(w^4,c_3^2,c_2c_3,c_2^3 w,c_2^4)$.
 A polynomial is said to be \emph{top-reducible} if any of its
 monomials is divisible by one of these leading terms; if so, we can
 subtract off a multiple of the corresponding relation to cancel the
 monomial, a process called \emph{top-reduction}.  Clearly, if a
 polynomial can be reduced to zero by iterated top-reduction then it
 must lie in $J$, but the converse need not hold for an arbitrary list
 of generators of an arbitrary ideal.  Let $a$ and $b$ be any two of
 our relations, let $a'$ and $b'$ be their leading terms, and let $c'$
 be the greatest common divisor of $a'$ and $b'$.  The corresponding
 \emph{syzygy} is the element $c:=(a'/c')b-(b'/c')a\in J$.  To say
 that our relations form a Gr\"obner basis means precisely that all
 these syzygies can be reduced to zero by iterated top-reduction.
 This can be checked by direct computation.  For example, the syzygy
 of $r'_4$ and $r'_2$ is the element
 $c_2r'_4-wr'_2=c_2^3w^3+c_2c_3w+c_3w^3$.  The first monomial is
 divisible by the leading term of $r'_4$, so we can top-reduce by
 subtracting $w^2 r'_4$ to get $c_2c_3w+c_2w^4$.  We can then do two
 more top-reductions by subtracting $w$ times the relation $c_2c_3$
 and $c_2$ times the relation $w^4$ to get $0$, as required.  Now
 observe that the $17$ monomials listed in the statement of the
 theorem are precisely those that are not top-reducible.  It follows
 from the theory of Gr\"obner bases 
 that they form a basis for $C(E,\Sg_4)$, as claimed.
\end{proof}

\begin{corollary}
 $C(K,\Sg_4)$ is a Gorenstein ring, and the element $w^3c_3$ generates
 the socle.
\end{corollary}
\begin{proof}
 One sees easily from the relations listed that $w$, $c_2$ and $c_3$
 annihilate $w^3c_3$, so $w^3c_3$ lies in the socle.  Now let $a$ be
 an arbitrary element of the socle.  It will be convenient to put
 $e=c_2^3+wc_2+c_3$ (so that $we=c_3e=0$) and to use the basis given
 in the Proposition but with $c_2^3$ replaced by $e$.  Using the
 equation $wa=0$ we see immediately that $a$ lies in the span of
 $\{w^3,w^3c_2,w^3c_2^2,e,w^3c_3\}$.  Using the equation $c_3a=0$ and
 the fact that $c_3c_2=c_3e=c_3^2=0$ we find that the coefficient of
 $w^3$ is zero, so $a=\al w^3c_2+\bt w^3c_2^2+\gm e+\dl w^3c_3$ say.
 One can check that $w^3c_2^3=w^3c_3$ and $c_2e=w^3c_2$, so 
 \[ 0 = c_2 a = \al w^3 c_2^2 + \bt w^3 c_3 + \gm w^3 c_2, \]
 so $\al=\bt=\gm=0$, so $a=\dl w^3c_3$.  This shows that the socle is
 one-dimensional, so the ring is Gorenstein as claimed.
\end{proof}

\subsection{A transfer argument}\label{subsec-transfer}

\begin{proposition}\label{prop-K-inj}
 The map $\tht\:C(K,\Sg_4)\xra{}K^0B\Sg_4$ is injective.
\end{proposition}
\begin{proof}
 Note that every nontrivial ideal in $C(K,\Sg_4)$ contains the socle,
 so it will suffice to show that the socle is not contained in
 $\ker(\tht)$, or equivalently that $w^3c_3\neq 0$ in $K^0B\Sg_4$.
 Let $P$ be the Sylow subgroup in $\Sg_4$; it will be enough to show
 that $w^3c_3$ has nontrivial image in $K^0BP$.  Put $V=P\cap A_4$;
 one can check that this consists of the identity and the three
 transposition pairs, so it is isomorphic to $C_2^2$.  Recall that the
 series $\dts(x)$ is defined to be $[2](x)/x$, which in our case is
 just $x^3$.  As $w$ is the Euler class of $\ep$ and
 $V=\ker(\ep\:P\xra{}C_2)$, standard arguments show that
 $\tr_V^P(1)=\dts(w)=w^3$.  This means that $w^3c_3=\tr_V^P(c_3)$.  To
 see that this is nonzero, we use the canonical bilinear form on
 $K^0BP$ defined in~\cite{st:kld}.  This satisfies Frobenius
 reciprocity, so $(\tr_V^P(c_3),1)_P=(c_3,1)_V$.  If we let $x$ and
 $y$ be the Euler classes of two of the nontrivial characters of $V$,
 then $K^0BV=\F_4[x,y]/(x^4,y^4)$ and the Euler class of the third
 character is $x+_Fy=x+y+x^2y^2$.  One checks that the restriction of
 $\rho$ to $V$ is the regular representation minus the trivial
 representation, which is the sum of the three nontrivial characters.
 This implies that the restriction of $c_3$ to $V$ is
 $xy(x+_Fy)=x^2y+xy^2+x^3y^3$.  Using~\cite[Corollary 9.3]{st:kld} we
 see that $(x^iy^j,1)_V$ is $1$ if $i=j=3$ and $0$ otherwise, so
 $(c_3,1)_V=1$.  As $(w^3c_3,1)_P=(c_3,1)_V=1$ we see that
 $w^3c_3\neq 0$, as claimed.
\end{proof}

\subsection{The proof of Theorem~\ref{thm-Sg-four}}
\label{subsec-proof} 

This is now easy.  We know from~\cite{st:msg} that $E^0B\Sg_4$ is a
free module of finite rank over $E^0$.  It follows by well-known
arguments that $K^0B\Sg_4=(E^0B\Sg_4)/I_2$, which is free of the same
rank over $K^0=E^0/I_2=\F_4$.  The rank is also the same as the rank
of $L\ot E^0B\Sg_4$ over $L$, and generalised character theory tells
us that this is equal to $|\Om(\Sg_4)|=17$.  Thus, the source and
target of the map $\tht\:C(K,\Sg_4)\xra{}K^0B\Sg_4$ both have rank
$17$ over $\F_4$ and Proposition~\ref{prop-K-inj} tells us that the
map is injective, so it must be an isomorphism.  Now consider the map
$\tht\:C(E,\Sg_4)\xra{}E^0B\Sg_4$.  This is an isomorphism modulo
$I_2$, so by Nakayama's lemma it is surjective.  As $E^0B\Sg_4$ is
free it is a split surjection, so $C(E,\Sg_4)=E^0B\Sg_4\oplus N$ say.
This implies that $C(K,\Sg_4)=K^0B\Sg_4\oplus N/I_2N$, so by counting
ranks we see that $N/I_2N=0$, so by Nakayama again we see that $N=0$.
Thus $C(E,\Sg_4)=E^0B\Sg_4$ as claimed.  We know from
Proposition~\ref{prop-K-Sg-four} that our list of $17$ monomials is a
basis for $K^0B\Sg_4$ over $\F_4$, and it now follows that it is also
a basis for $E^0B\Sg_4$ over $E^0$.

\section{Extraspecial $p$-groups}\label{sec-xspec}

In this section we define a class of ``extraspecial'' $p$-groups
(where $p$ is an odd prime), and show that for these groups the map
$\kp\:\Om(G)\xra{}\OmCh(G)$ is injective but not surjective.  It
follows using Theorem~\ref{thm-gen-char} that the map
$\tht\:C(E,G)\xra{}E^0BG$ cannot be an isomorphism.  We have not
investigated the situation more deeply than this.

Let $V$ be an elementary Abelian $p$-group of rank $2d$ equipped with
a nondegenerate alternating form $b\:V\tm V\xra{}\F_p$.  We will say
that a subspace $W\leq V$ is \emph{isotropic} if $b(u,v)=0$ for all
$u,v\in W$.

Let $G$ be the set $\F_p\tm V$ with the group operation
$(x,u).(y,v)=(x+y+b(u,v),u+v)$.  This has order $p^{2d+1}$ and
exponent $p$, and it fits in a central extension
\[ Z=\F_p \xra{j} G \xra{q} V. \]
In fact $Z$ is the centre of $G$, and the non-central conjugacy
classes are the fibres of $q$ over $V\sm\{0\}$, so they all have order
$p$.  This gives $p+p^{2d}-1$ conjugacy classes altogether.

We can evidently view $R(V)=\Z[V^*]$ as a sub $\Lm$-ring of $R(G)$.

\begin{definition}
 For any nontrivial character $\zt\:Z\xra{}S^1$, let $\phi(\zt)$ be
 the class function on $G$ defined by
 \[ \phi(\zt)(g) = \begin{cases}
     p^d\zt(g) & \text{ if } g\in Z \\
     0         & \text{ otherwise. }
     \end{cases}
 \]
 We also write $\rho_V$ for the regular representation of $V$, and
 $\rho_G$ for the regular representation of $G$.
\end{definition}

The following result is standard, but we give a proof for
completeness. 
\begin{proposition}
 For each $\zt\in Z^*\sm\{1\}$, the class function $\phi(\zt)$ is an
 irreducible character.  Moreover, we have
 \[ R(G) = \Z[V^*]\oplus\Z\{\phi(\zt)\st\zt\in Z^*\sm\{1\}\}. \]
\end{proposition}
\begin{proof}
 Choose a maximal isotropic subspace $W\leq V$, so $W\simeq\F_p^d$.
 Put $H=q^{-1}W\leq G$, which is isomorphic to $Z\tm W$ as a group
 because $W$ is isotropic.  Let $r\:H\xra{}Z$ be the projection and
 put $\sg=\ind_H^G\pi^*\zt$.  We claim that $\sg=\phi(\zt)$.  To see
 this, first note that $H$ is normal in $G$, so $\sg(g)=0$ for
 $g\not\in H$.  Next, suppose that $g\in H\sm Z$, say $g=(x,w)$ with
 $w\in W\sm\{0\}$.  Let $U$ be such that $V=W\oplus U$, so
 $U\simeq\F_p^d$ and $(0,u)^{-1}(x,w)(0,u)=(x-2b(u,w),w)$.  From the
 definitions we see that $\sg(x,w)=\sum_u\zt(x-2b(u,w),0)$.  The map
 $u\mapsto (x-2b(u,w),0)$ is a surjection from $U$ to $Z$, each of
 whose fibres has the same order, and $\zt\:Z\xra{}S^1$ is a
 nontrivial homomorphism; it follows easily that $\sg(x,w)=0$, as
 required.  Finally, suppose that $g\in Z$, say $g=(x,0)$.  Then
 $(0,u)^{-1}(x,0)(0,u)=(x,0)$ so
 $\sg(x,0)=\sum_u\zt(x,0)=p^d\zt(x,0)$.  This shows that
 $\sg=\phi(\zt)$ as claimed, so $\phi(\zt)$ is a character.  One
 checks easily that
 $\langle\phi(\zt),\phi(\zt)\rangle=|G|^{-1}\sum_{z\in Z}p^{2d}=1$, so
 $\phi(\zt)$ is irreducible.  As $\zt$ runs over $Z^*\sm\{1\}$ this
 gives $p-1$ distinct irreducibles of degree $p^d$, and $V^*$ gives a
 further $p^{2d}$ distinct irreducibles of degree $1$.  We have seen
 that $G$ has $p^{2d}+p-1$ conjugacy classes and thus $p^{2d}+p-1$
 irreducible characters, so our list is complete.  It follows that 
 $R(G)=\Z[V^*]\oplus\Z\{\phi(\zt)\st\zt\in Z^*\sm\{1\}\}$ as claimed.
\end{proof}

\begin{lemma}\label{lem-lm-rho}
 Let $C$ be cyclic of order $p$.  Then 
 \[ \lm^k(p^{d-1}\rho_C) = \begin{cases}
    \bsm p^{d-1} \\ k/p \esm + 
    \frac{1}{p}
      \left(\bsm p^d\\ k\esm - \bsm p^{d-1}\\ k/p \esm \right) \rho_C
    & \text{ if } p|k \\
    \frac{1}{p} \bsm p^d \\ k \esm \rho_C
    & \text{ otherwise }
   \end{cases}
 \]     
\end{lemma}
\begin{proof}
 Let $\chi$ be a generator of $C^*$, so $R(C)=\Z[\chi]/(\chi^p-1)$ and
 $\rho_C=\sum_{j=0}^{p-1}\chi^j$.  We have $\chi\rho_C=\rho_C$ and so
 $\lm^k(\rho_C)=\lm^k(\chi\rho_C)=\chi^k\lm^k(\rho_C)$.  If $0<k<p$
 then $\chi^k$ is also a generator, and it follows that
 $\lm^k(\rho_C)$ is an integer multiple of $\rho_C$.  On the other
 hand, it is easy to check that $\lm^0(\rho_C)=\lm^p(\rho_C)=1$.  If
 we put $A=\Z\{1,\rho_C\}$ then $A$ is a subring of $R(C)$ (with
 $\rho_C^2=p\rho_C$) and $\lm_t(\rho_C)\in A[t]$ so
 $\lm_t(p^{d-1}\rho_C)=\lm_t(\rho_C)^{p^{d-1}}$ also lies in $A[t]$,
 say $\lm^k(\rho_C)=n_k+m_k\rho_C$.  Moreover, if we work mod $\rho_C$
 we have $\lm_t(\rho_C)\cong 1+t^p$ so
 $\lm_t(p^{d-1}\rho_C)\cong(1+t^p)^{p^{d-1}}$.  Thus, if $p$ divides
 $k$ then $n_k=\bcf{p^{d-1}}{k/p}$, and if $p$ does not divide $k$
 then $n_k=0$.  Moreover, by counting dimensions we see that
 $n_k+pm_k=\bcf{p^d}{k}$ for all $k$.  The lemma now follows easily.
\end{proof}

\begin{proposition}
 \begin{align*}
  \chi\phi(\zt) &= \phi(\zt) \\
  \phi(\zt)\phi(\xi) &= 
   \begin{cases} \rho_V & \text{ if } \zt\xi=1 \\
                 p^d\phi(\zt\xi) & \text{ otherwise }
   \end{cases} \\
  \psi^k\chi &= \chi^k \\
  \psi^k\phi(\zt) &= 
   \begin{cases} p^d & \text{ if } p|k \\
                 \phi(\zt^k) & \text{ otherwise }
   \end{cases} \\
  \lm^k(\phi(\zt)) &= 
   \begin{cases}
    \bsm p^{d-1} \\ k/p \esm + 
    \frac{1}{p^{2d}}
      \left(\bsm p^d\\ k\esm - \bsm p^{d-1}\\ k/p \esm \right) \rho_V
    & \text{ if } p|k \\
    \frac{1}{p^d} \bsm p^d \\ k \esm \phi(\zt^k) 
    & \text{ otherwise }
   \end{cases} \\
 \end{align*}
\end{proposition}
\begin{proof}
 Everything except for $\lm^k(\phi(\zt))$ can be done by easy
 manipulation of characters.  For the remaining case, it suffices to
 check that the claimed equations hold when restricted to any cyclic
 subgroup $C\leq G$.  First consider the case $C=Z$, so $\rho_V$
 restricts on $C$ to the trivial representation of degree $p^{2d}$.
 Then $\phi(\zt)$ becomes $p^d\zt$, so $\lm^k\phi(\zt)$ becomes
 $\bcf{p^d}{k}\zt^k$.  Using this, it is easy to check that the
 equations hold when restricted to $Z$.

 Now suppose instead that $C\leq G$ is a cyclic group not contained in
 $Z$ (which implies that $|C|=p$).  Then $\rho_V|_C=p^{2d-1}\rho_C$
 and $\phi(\xi)|_C=p^{d-1}\rho_C$ for all $\xi\in Z^*\sm\{1\}$.  
 Using Lemma~\ref{lem-lm-rho} we deduce that our equations for
 $\lm^k(\phi(\zt))$ are correct when restricted to $C$, as required.
\end{proof}

\begin{definition}
 For any homomorphism $\al\:V^*\xra{}\Tht$, put
 \[ c_\al=\sum_{\chi\in V^*}[\al(\chi)]\in\Z[\Tht(1)]_{p^{2d}}^+ \]
 and
 \[ U_\al = \{u\in\Z[\Tht(1)]_{p^d}^+\st u\psi^{p-1}(u)=c_\al\}.\]
 We also put $U=\{(\al,u)\st u\in U_\al\}$.
\end{definition}

\begin{theorem}\label{thm-OmCh-U}
 There is a natural bijection $\OmCh(G)=U$.  The map
 $\kp\:\Om(G)\xra{}U$ is injective, and the image is the set of pairs
 $(\al,u)\in U$ such that the image of the dual map
 $\al^*\:\Tht^*\xra{}V$ is isotropic.
\end{theorem}

The proof will follow after a lemma.

\begin{lemma}\label{lem-U}
 Let $\al\:V^*\xra{}\Tht$ be a homomorphism with image $A$ of order
 $p^e$.  If $e>d$ then $U_\al=\emptyset$.  If $e\leq d$ then $U_\al$
 is the set of elements of the form $u=p^{d-e}\sum_{c\in C}[c]$, where
 $C$ runs over the cosets of $A$ in $\Tht(1)$.
\end{lemma}
\begin{proof}
 Put $c'_\al=\sum_{a\in A}[a]\in\Z[\Tht(1)]_{p^e}^+$ so that
 $c_\al=p^{2d-e}c'_\al$.  Suppose $u\in U_\al$ and that $b\in u$.  Put
 $v=[-b]u$, so $v\in U_\al$ and $0\in v$.  Thus $0\in\psi^{p-1}(v)$
 also, so $v\leq v\psi^{p-1}(v)=c_\al=p^{2d-e}c'_\al$, so we can write
 $v=\sum_{a\in A}n_a[a]$ for suitable natural numbers $n_a$.  By
 looking at the multiplicity of $[0]$ in the equation
 $v\psi^{p-1}(v)=p^{2d-e}c'_\al$ we see that $\sum_an_a^2=p^{2d-e}$.
 On the other hand, as $v\in\Z[\Tht(1)]_{p^d}^+$ we have
 $\sum_an_a=p^d$.  It follows that
 \[ \sum_a(n_a-p^{d-e})^2=
    \sum_an_a^2-2p^{d-e}\sum_an_a+p^{2d-2e}\sum_a1=
    p^{2d-e}-2p^{2d-e}+p^{2d-e}=0,
 \]
 so $n_a=p^{d-e}$ for all $a$.  If we now let $C$ be the coset $b+A$
 we find that $u=p^{d-e}\sum_{c\in C}[c]$.  Conversely, it is trivial
 to check that any element of this form lies in $U_\al$.
\end{proof}

\begin{proof}[Proof of Theorem~\ref{thm-OmCh-U}]
 Let $\zt$ be the usual character $x\mapsto e^{2\pi ix/p}$ of
 $Z=\Z/p$.  Given $f\:R(G)\xra{}\Z[\Tht]$ it is clear that the
 restriction of $f$ to $R(V)=\Z[V^*]$ gives a homomorphism
 $\al\:V^*\xra{}\Tht(1)\leq\Tht$, and we put
 $u=f(\phi(\zt))\in\Z[\Tht]_{p^d}^+$.  As $f$ is a $\Lm$-ring
 homomorphism we have
 \[ u\psi^{p-1}(u) = 
    f(\phi(\zt)\phi(\zt^{p-1})) = 
    f(\rho_V) =
    \sum_\chi [\al(\chi)],
 \]
 so $(\al,u)\in U$.

 Conversely, suppose we start with $(\al,u)\in U$.  Let $e$ and $C$ be
 as in Lemma~\ref{lem-U}.  We define a homomorphism
 $f\:R(G)\xra{}\Z[\Tht]$ of additive groups by $f(\chi)=[\al(\chi)]$
 for $\chi\in V^*$ and 
 \[ f(\phi(\zt^k)) = \psi^k(u) = p^{d-e}\sum_{c\in kC}[c] \]
 for $k\in\Z\sm p\Z$.  It is easy to check that this is a ring
 homomorphism that sends $R^+_k(G)$ to $\Z[\Tht]^+_k$ and commutes
 with the Adams operations.  As $\Z[\Tht]$ is torsion free it follows
 that $f$ commutes with $\Lm$-operations as well, so $f\in\OmCh(G)$.
 Clearly these constructions give the required bijection $\OmCh(G)=U$.

 Now suppose we have a homomorphism $\mu\:\Tht^*\xra{}G$.  Then
 $\mu(u)=(\om(u),\sg(u))$ for some functions $\om\:\Tht^*\xra{}\F_p$
 and $\sg\:\Tht^*\xra{}V$.  As $\mu$ and the projection $q\:G\xra{}V$
 are homomorphisms we see tht $\sg$ is a homomorphism.  Let $W\leq V$
 be the image of $\sg$, and put $e=\dim_{\F_p}W$.  As the image of
 $\mu$ must be commutative, it is not hard to see that $W$ is
 isotropic, so $e\leq d$.  As $q^{-1}W\simeq \F_p\tm W$ as groups, we
 see that $\om$ is also a homomorphism.  If we conjugate $(\om,\sg)$
 by $(x,u)\in G$ we get the homomorphism $(\om+\tau,\sg)$ where
 $\tau(t)=2b(u,\sg(t))$.  As $b$ is a perfect pairing, $\tau$ can be
 any map $\Tht^*\xra{}\F_p$ that factors through $\sg$, so $(\om,\sg)$
 is conjugate to $(\om',\sg)$ if and only if
 $\om|_{\ker(\sg)}=\om'|_{\ker(\sg)}$.  Now let $\sg^*\:V^*\xra{}\Tht$
 be the dual of $\sg$ and put $A=\sg^*(V^*)$, so $|A|=p^e$.  We also
 have a map $\om^*\:\F_p^*\xra{}\Tht$ and thus a point
 $t=\om^*(\zt)\in\Tht(1)$.  In $R(\F_p\tm W)=\Z[\F_p^*]\ot\Z[W^*]$ we
 have 
 \[ \phi(\zt)|_{\F_p\tm W}=p^{d-e}\zt\ot\rho_W=
     p^{d-e}\sum_{\xi\in W^*}\zt\ot\xi,
 \]
 and it follows that
 $\mu^*\phi(\zt)=p^{d-e}\sum_{a\in A}[t+a]\in\Z[\Tht]$.  Thus, if we
 write $[\mu]$ for the conjugacy class of $\mu$ then
 $\kp[\mu]=(\sg^*,p^{d-e}\sum_{a\in A}[t+a])\in U$.  It follows that
 $\kp[\mu]$ determines $\sg$, and it also determines $t$ modulo $A$,
 so it determines $\om$ modulo $\sg^*(\Hom(V,\F_p))$, so it determines
 the conjugacy class $[\mu]$.  This proves that $\kp$ is injective as
 claimed.  We leave it to the reader to check that the image is as
 described.  
\end{proof}

\section{An apparently more precise approach}\label{sec-precise}

There are some senses in which the $\Lm$-operations do not capture all
possible information about the representation theory of $G$, and it is
reasonable to wonder whether a more accurate approximation to $X(G)$
could be defined by taking more information into account.  In this
section we show that this is not the case: we construct an
approximation $Y(G)$ using all possible operations, and show that it
is the same as $\XCh(G)$.

\begin{definition}
 Let $\CG$ be the category of Lie groups and continuous homomorphisms,
 and let $\CGb$ be the quotient category in which conjugate
 homomorphisms are identified.  Let $\CN$ be the set of finite
 sequences $\unn=(n_1,\ldots,n_r)$ with $n_i\in\N$.  For $\unn\in\CN$
 we put $\GL(\unn)=\prod_i\GL(n_i,\C)$ and 
 \[ R(\unn,G)=\prod_iR^+_{n_i}(G)=\CGb(G,\GL(\unn)).  \]
 We make $\CN$ into a category by putting
 $\CN(\unn,\unm)=\CG(\GL(\unn),\GL(\unm))$, and we let $\CNb$ be the
 category with the same objects and with morphisms
 $\CNb(\unn,\unm)=\CGb(\GL(\unn),\GL(\unm))$; clearly this gives a
 covariant functor $\unn\mapsto R(\unn,G)$ from $\CNb$ to sets.

 Next, let $T(\unn)$ be the evident maximal compact torus in
 $\GL(\unn)$, so $T(\unn)\simeq\prod_{j=1}^NS^1$ where
 $N=\sum_{i=1}^rn_i$.  Let $W(\unn)$ be the Weyl group of $T(\unn)$,
 so $W(\unn)\simeq\prod_i\Sg_{n_i}$.  We can thus form the scheme
 \[ D(\unn)=\Hom(T(\unn)^*,\GG)/W(\unn)=\prod_i\Div_{n_i}^+(\GG). \]
 By elementary arguments in representation theory we see that any
 homomorphism $f\:\GL(\unn)\xra{}\GL(\unm)$ is conjugate to one that
 sends $T(\unn)$ into $T(\unm)$, and that the resulting map
 $T(\unn)\xra{}T(\unm)$ is unique up to the action of $W(\unm)$.  
 Using this, it is easy to make the assignment $\unn\mapsto D(\unn)$
 into a functor $\CNb\xra{}\hX_X$.

 Finally, we define a functor $Y(G)$ from discrete $\OX$-algebras to
 sets by putting
 \[ Y(G)(A) = [\CNb,\text{Sets}](R(-,G),D(-)(A)). \]
\end{definition}

\begin{theorem}\label{thm-YG}
 There is a natural isomorphism $Y(G)\simeq\XCh(G)$.
\end{theorem}

Before proving this, we relate $Y(G)$ to an auxiliary model involving
unitary groups rather than general linear groups.
\begin{definition}
 Let $\CGt$ be the quotient of $\CG$ in which homomorphisms
 $u,v\:U\xra{}V$ are identified if $u|_K$ is conjugate to $v|_K$ for
 every compact subgroup $K\leq U$.  (For example, the homomorphism
 $u\:\GL(1)\xra{}\GL(1)$ given by $u(z)=|z|$ becomes trivial in
 $\CGt$.)  Let $\CNt$ be the category with the same objects as $\CN$
 and morphisms $\CNt(\unn,\unm)=\CGt(\GL(\unn),\GL(\unm))$.  As $G$
 and $T(\unn)$ are compact, it is clear that the functors $R(-,G)$ and
 $D(-)$ factor through $\CNt$, and thus that
 \[ Y(G)(A) = [\CNt,\text{Sets}](R(-,G),D(-)(A)). \]
\end{definition}

\begin{lemma}\label{lem-conjugate-in-U}
 Let $K$ be a compact Lie group, and let $v,w\:K\xra{}U(d)$ be
 continuous homomorphisms.  If $v$ and $w$ are conjugate in $\GL(d)$,
 then they are conjugate in $U(d)$.
\end{lemma}
\begin{proof}
 The statement can easily be translated as follows: Let $V$ and $W$ be
 finite-dimensional vector spaces over $\C$ equipped with actions of
 $G$ and invariant Hermitian inner products.  Then if there exists an
 equivariant isomorphism $f\:V\xra{}W$, then $f$ can be chosen to
 preserve the inner products.

 To see this, we first recall some facts about invariant Hermitian
 products.  For any complex vector space $V$ we let $\bV$ be the same
 set with the conjugate action of $\C$, and let $\bV^*$ be the dual of
 $\bV$.  The set of Hermitian products $\bt$ on $V$ bijects with the
 set of isomorphisms $\bt'\:V\xra{}\bV^*$ satisfying certain symmetry
 and positivity conditions.  For any representation $V$ one can always
 choose an invariant Hermitian product so $V$ is equivariantly
 isomorphic to $\bV^*$.  For each irreducible representation $S$ we
 fix a Hermitian product $\bt_S$ on $S$; Schur's lemma implies that
 $\Hom_K(S,\ov{S}^*)=\C\bt'_S$ and that any other invariant Hermitian
 product is a positive scalar multiple of $\bt_S$.

 Now let $\bt$ be a Hermitian product on $V$ and suppose that
 $V=V_0\oplus V_1$ and $\Hom_K(V_1,V_0)=0$.  Then
 $\Hom_K(V_1,\bV_0^*)=0$ and
 \[ \Hom_K(V_0,\bV_1^*)=\Hom_K(\bV_1^*,V_0)^*=\Hom_K(V_1,V_0)^*=0 \]
 so the equivariant isomorphism $\bt'\:V\xra{}\bV^*$ must have the
 form $\bt'_0\oplus\bt'_1$ for some $\bt'_i\:V_i\xra{}\bV_i^*$.  This
 implies that $V_0$ and $V_1$ are orthogonal with respect to $\bt$.

 Now let $S_1,\ldots,S_t$ be the distinct irreducible representations
 that occur in $V$.  Then there is a unique decomposition
 $V=V_1\oplus\ldots\oplus V_t$, where $V_i\simeq\C^{d_i}\ot S_i$ for
 some $d_i$ and thus $\Hom_K(V_i,V_j)=0$ when $i\neq j$.  By the
 previous paragraph, the subspaces $V_i$ are orthogonal to each other.
 As $\Hom_K(S_i,\ov{S}_i^*)=\C\bt_{S_i}$, we find that the restriction
 of $\bt$ to $V_i$ has the form $\gm_i\ot\bt_{S_i}$ for some Hermitian
 product $\gm_i$ on $\C^{d_i}$.  By Gram-Schmidt, the space
 $(\C^{d_i},\gm_i)$ is isomorphic to $\C^{d_i}$ with its usual
 Hermitian product, so $(V_i,\bt|_{V_i})$ is equivariantly and
 isometrically isomorphic to the orthogonal direct sum of $d_i$ copies
 of $(S_i,\bt_{S_i})$.  This means that the numbers $d_i$ determine
 the \emph{isometric} isomorphism type of $V$, and the lemma follows
 immediately.
\end{proof}

\begin{lemma}\label{lem-unitary}
 There are natural bijections
 \[ \CNt(\unn,\unm)=\CGb(U(\unn),\GL(\unm))=\CGb(U(\unn),U(\unm)),
 \]
 where $U(\unn)=\prod_iU(n_i)\leq\GL(\unn)$.
\end{lemma}
\begin{proof}
 It is easy to reduce to the case where the list $\unm$ has length
 $1$, say $\unm=(d)$.  As any representation of $U(\unn)$ admits a
 Hermitian inner product, we see that the map
 $\CGb(U(\unn),U(d))\xra{}\CGb(U(\unn),\GL(d))$ is surjective.  It is
 also injective by Lemma~\ref{lem-conjugate-in-U}.  Similarly, by
 considering invariant Hermitian products we see that if $K$ is
 compact and $u\:K\xra{}\GL(\unn)$ then $u$ is conjugate to a
 homomorphism $K\xra{}U(\unn)$.  By applying this to the inclusion
 map, we see that any compact subgroup of $\GL(\unn)$ is conjugate to
 a subgroup of $U(\unn)$.  It follows that any two homomorphisms
 $v,w\:\GL(\unn)\xra{}\GL(d)$ are identified in
 $\CGt(\GL(\unn),\GL(d))$ iff their restrictions to $U(\unn)$ are
 conjugate, so we have a well-defined and injective restriction map
 $\CNt(\unn,d)\xra{}\CGb(U(\unn),\GL(d))$.  It is an easy consequence
 of the theory of roots and so on that any representation of $U(\unn)$
 extends uniquely to a complex-analytic representation of $\GL(\unn)$,
 so our restriction map is also surjective.
\end{proof}

\begin{proof}[Proof of Theorem~\ref{thm-YG}]
 Consider a point $g\in Y(G)(A)$, in other words a natural
 transformation $g_\unn\:R(\unn,G)\xra{}D(\unn)(A)$ for $\unn\in\CN$.
 By putting together the maps
 \[ g_d\:R^+_d(G)=R(d,G)\xra{}D(d)(A)=\Div_d^+(\GG)(A), \]
 we get a function $f\:R^+(G)\xra{}\Div^+(\GG)(A)$.  Next, for any
 $d,e\geq 0$ we have projections $\GL(d)\xla{}\GL(d,e)\xra{}\GL(e)$
 and we can use the resulting maps to identify $R((d,e),G)$ with
 $R_d^+(G)\tm R_e^+(G)$ and $D(d,e)(A)$ with
 $\Div_d^+(\GG)(A)\tm\Div_e^+(\GG)(A)$ and $g_{(d,e)}$ with
 $g_d\tm g_e$.  There are evident maps
 $\oplus\:\GL(d,e)\xra{}\GL(d+e)$ and $\ot\:\GL(d,e)\xra{}\GL(de)$,
 and using the naturality of $g$ with respect to these maps we find
 that $f$ is a semiring homomorphism.  Similarly, we have maps
 $\lm^k\:\GL(d)\xra{}\GL(\bcf{d}{k})$ in $\CGb$ and the naturality of
 $g$ with respect to these maps implies that $f$ commutes with
 $\Lm$-operations.  It is clear that
 $f(R_d^+(G))\sse\Div_d^+(\GG)(A)$, so $f\in\XCh(G)(A)$.  We define a
 map $\rho\:Y(G)\xra{}\XCh(G)$ by $\rho(g)=f$.  Because
 $g_\unn=g_{n_1}\tm\ldots\tm g_{n_r}$ we see that $\rho$ is
 injective. 

 Now suppose we start with a point $f\in\XCh(G)(A)$.  Let
 $g_d\:R(d,G)\xra{}D(d)(A)$ be the restriction of
 $f\:R(G)\xra{}\Div(\GG)(A)$, and put 
 \[ g_\unn=g_{n_1}\tm\ldots g_{n_r} \: R(\unn,G)\xra{}D(\unn)(A). \]
 We need to check that this gives a natural transformation.  As
 $R(\unm,G)=\prod_iR(m_i,G)$ and $D(\unm)(A)=\prod_iD(m_i)(A)$, it
 suffices to check naturality for maps $u\:\unn\xra{}d$ in $\CNt$, or
 equivalently (by Lemma~\ref{lem-unitary}) for homomorphisms
 $u\:U(\unn)\xra{}\GL(d)$.  We need to show that the left hand square
 in the following diagram commutes:
 \begin{diag}
  \node{R(\unn,G)} 
  \arrow{e,t}{u_*}
  \arrow{s,l}{g_{\unn}}
  \node{R_d^+(G)}
  \arrow{s,l}{g_d}
  \arrow{e,V}
  \node{R(G)}
  \arrow{s,r}{f} \\
  \node{D(\unn)(A)}
  \arrow{e,b}{u_*}
  \node{\Div_d^+(\GG)}
  \arrow{e,V}
  \node{\Div(\GG)(A).}
 \end{diag}
 The right hand square commutes and the two right hand horizontal maps
 are injective so it suffices to show that the two composite maps
 $R(\unn,G)\xra{}\Div(\GG)(A)$ are the same.  We call these two maps
 $\al(u)$ and $\bt(u)$.  Let $F$ be the set of all functions from
 $R(\unn,G)$ to $\Div(\GG)(A)$, thought of as a ring with pointwise
 operations.  It is formal to check that $\al(u+v)=\al(u)+\al(v)$ and
 $\al(uv)=\al(u)\al(v)$, so $\al$ is a homomorphism of semirings from
 $R^+(U(\unn))$ to $\Div(\GG)(A)$.  It can thus be extended to a ring
 map $R(U(\unn))\xra{}\Div(\GG)(A)$, and the same applies to $\bt$.
 It is well-known that $R(U(\unn))=\bigotimes_iR(U(n_i))$ so it
 suffices to check that $\al=\bt$ on $R(U(n_i))$ for all $i$.  This
 reduces us to the case where $\unn=(e)$ say.  It is also well-known
 that $R(U(e))=\Z[\lm^1,\ldots,\lm^e][(\lm^e)^{-1}]$, so it suffices
 to check that $\al(\lm^j)=\bt(\lm^j)$, which is true because $f$ is a
 homomorphism of $\Lm$-rings.  

 This shows that $g\in Y(G)(A)$, and clearly $\rho(g)=f$.  Thus $\rho$
 is surjective and hence an isomorphism.
\end{proof}

\section{A result on restrictions of characters}\label{sec-restriction}

\begin{theorem}\label{thm-restriction}
 Let $G$ be a finite group with a normal subgroup $N$ such that $|N|$
 is coprime to $|G/N|$.  Then the restriction map
 $R^+(G)\xra{}R^+(N)^G$ is surjective.
\end{theorem}
The proof will follow after some preliminary results.

\begin{lemma}\label{lem-lift}
 Let $H$ be a group, and let $W,X,Y$ be $H$-sets, with equivariant
 maps $W\xra{f}X\xla{q}Y$.  Then there is an equivariant map
 $\tilde{f}\:W\xra{}Y$ with $q\tilde{f}=f$ iff for each $w\in W$ there
 exists $y\in Y$ with $q(y)=f(w)$ and $\stab_H(y)\geq\stab_H(w)$.
\end{lemma}
\begin{proof}
 Write $W$ as a disjoint union of orbits.
\end{proof}

\begin{lemma}\label{lem-restriction}
 Let $G$ and $N$ be as above, and let $\rho\:N\xra{}GL(V)$ be an
 irreducible representation of $N$ whose character is stable under
 $G$.  Then there is a homomorphism $\sg\:G\xra{}GL(V)$ extending
 $\rho$.  
\end{lemma}
\begin{proof}
 Suppose $g\in G$, and define $\rho^g\:N\xra{}GL(V)$ by
 $\rho^g(x)=\rho(gxg^{-1})$.  By hypothesis, this has the same
 character as $\rho$, so there exists an intertwining operator
 $\tht\:V\xra{}V$ such that $\rho^g(x)=\tht^{-1}\rho(x)\tht$ for all
 $x\in N$.  As $V$ is an irreducible representation of $N$ we see that
 $\Aut_N(V)=\C$ and thus $\tht$ is unique up to multiplication by a
 scalar matrix.  We can thus define a map $\phi\:G\xra{}PGL(V)$ by
 $\phi(g)=[\tht]$; this is a homomorphism making the following diagram
 commute. 
 \begin{diag}
  \node{N} \arrow{s,l}{\rho} \arrow{e,V} 
  \node{G} \arrow{s,r}{\phi} \\
  \node{GL(V)} \arrow{e,b,A}{\pi} \node{PGL(V).}
 \end{diag}
 Put $n=|N|$ and $d=\dim_\C(V)$.  As $V$ is irreducible we know that
 $d$ divides $n$.  Put $Y=\{\al\in GL(V)\st\det(\al)^n=1\}$, and note
 that $\pi\:Y\xra{}PGL(V)$ is surjective and $\rho(N)\leq Y$.  Let
 $N^2$ act on $G$ by $(x,y).g=xgy^{-1}$ and on $GL(V)$ by
 $(x,y).\al=\rho(x)\al\rho(y)^{-1}$.  

 We claim that there is an $N^2$-equivariant map $\zt\:G\xra{}Y$ such
 that $\pi\zt=\phi$ and $\zt=\rho$ on $N$.  Clearly
 $G=N\amalg(G\sm N)$ as $N^2$-sets and $\rho\:N\xra{}Y$ is
 $N^2$-equivariant, so it suffices to define $\zt$ on $G\sm N$.  Fix
 $g\in G\sm N$, and choose $\tht$ as before.  After multiplying by a
 suitable scalar, we may assume that $\det(\tht)=1$ so $\tht\in Y$.
 By Lemma~\ref{lem-lift}, it will suffice to show that
 $\stab_{N^2}(g)\leq\stab_{N^2}(\tht)$.  Suppose that $(x,y)$
 stabilises $g$, so $xgy^{-1}=g$, so $y=g^{-1}xg$.  By the definition
 of $\tht$ we have $\rho(y)=\tht^{-1}\rho(x)\tht$, or in other words
 $\rho(x)\tht\rho(y)^{-1}=\tht$, so $(x,y)$ stabilises $\tht$, as
 required.

 Now define $\xi\:G^2\xra{}Y$ by $\xi(g,h)=\zt(h)\zt(gh)^{-1}\zt(g)$.
 Clearly $\pi\xi(g,h)=1$, and the kernel of $\pi\:Y\xra{}GL(V)$ is the
 group $C_{nd}$ of $nd$'th roots of unity, so we can regard $\xi$ as a
 map $G^2\xra{}C_{nd}$.  As $\zt$ is equivariant, it is easy to check
 that $\xi(xg,hy)=\xi(g,h)$ for $x,y\in N$, so we get an induced map
 $\xib\:(G/N)^2\xra{}C_{nd}$.  One also sees directly that for
 $g,h,k\in G/N$ we have
 \[ \xib(h,k)\xib(gh,k)^{-1}\xib(g,hk)\xib(g,h)^{-1}=1, \]
 so $\xib$ is a $2$-cocycle.  On the other hand $nd$ divides $n^2$ and
 thus is coprime to $|G/N|$, so we have $H^2(G/N;C_{nd})=0$.  We can
 thus choose a function $\om\:G/N\xra{}C_{nd}$ such that
 $\xi(g,h)=\om(h)\om(gh)^{-1}\om(g)$ for all $g,h\in G$.  By putting
 $g=h=1$ we see that $\om(1)=1$ and thus $\om(x)=1$ for $x\in N$.  We
 define $\sg(g)=\om(g)^{-1}\zt(g)$; this clearly gives a homomorphism
 $G\xra{}GL(V)$ with $\sg|_N=\rho$, as required.
\end{proof}

\begin{proof}[Proof of Theorem~\ref{thm-restriction}]
 For each irreducible representation $\rho$ of $N$, let $\rho'$ denote
 the sum of the inequivalent $G$-conjugates of $\rho$.  Any
 $G$-invariant character is a direct sum of copies of the characters
 of the representations $\rho'$, so it suffices to show that $\rho'$
 extends to a representation of $G$.  Let $H$ be the stabiliser of
 $\chi_\rho$, so $N\leq H\leq G$.  Lemma~\ref{lem-restriction} implies
 that $\rho$ can be extended to a representation $\sg$ of $H$, and one
 sees from the Mackey formula that $\res^G_N\ind_H^G(\sg)\simeq\rho'$,
 so $\ind_H^G(\sg)$ is the required extension of $\rho'$.
\end{proof}


\end{document}